\documentclass[11pt]{article}
\usepackage{subeqn}
\usepackage{amssym}

\newcommand{\fbnp}{{f^b_n}\,'}

\newcommand{\X}{X}
\newcommand{\Y}{Y}
\newcommand{\x}{x}
\newcommand{\Q}{{\Bbb Q}}

\newcommand{\SSS}{{\mathbf S}}
\newcommand{\nichts}{{\left.\right.}}

\newcommand{\inv}{{\rm inv}}
\newcommand{\sech}{{\rm sech}}
\newcommand{\vp}{\epsilon}
\newcommand{\nn}{\nonumber}

\def\bj{{\bf j}}
\newcommand{\bs}{{\bf s}}
\def\ss{s}

\newtheorem{theorem}{Theorem}[section]
\newtheorem{lemma}[theorem]{Lemma}

\newtheorem{prop}[theorem]{Proposition}
\newtheorem{theo}[theorem]{Theorem}

\begin{document}

\newcommand{\ls}[1]
   {\dimen0=\fontdimen6\the\font \lineskip=#1\dimen0
\advance\lineskip.5\fontdimen5\the\font \advance\lineskip-\dimen0
\lineskiplimit=.9\lineskip \baselineskip=\lineskip
\advance\baselineskip\dimen0 \normallineskip\lineskip
\normallineskiplimit\lineskiplimit \normalbaselineskip\baselineskip
\ignorespaces }

\newcommand{\lbl}{\label}
\newcommand{\ignore}[1]{}{}
\def\squarebox#1{\hbox to #1{\hfill\vbox to #1{\vfill}}}
\newcommand{\qed}{\hspace*{\fill}
            \vbox{\hrule\hbox{\vrule\squarebox{.667em}\vrule}\hrule}\smallskip}
\newcommand{\oo}{\overline}
\newcommand{\uu}{\underline}

\newcommand{\be}{\begin{equation}}
\newcommand{\ee}{\end{equation}}
\newcommand{\bea}{\begin{eqnarray}}
\newcommand{\eea}{\end{eqnarray}}

\newcommand{\beq}[1]{\begin{equation}\label{#1}}
\newcommand{\eeq}{\end{equation}}
\newcommand{\beqn}[1]{\begin{eqnarray}\label{#1}}
\newcommand{\eeqn}{\end{eqnarray}}

\newcommand{\beaa}{\begin{eqnarray*}}
\newcommand{\eeaa}{\end{eqnarray*}}

\newcommand{\req}[1]{(\ref{#1})}

\newcommand{\proof}{\noindent{\bf Proof:}\ }
\newcommand{\al}{\alpha}
\newcommand{\ep}{\varepsilon}
\newcommand{\lam}{\lambda}
\newcommand{\la}{\lambda}
\newcommand{\La}{\Lambda}
\newcommand{\om}{\omega}
\newcommand{\reals}{{\Bbb R}}
\newcommand{\realsd}{\reals^d}
\newcommand{\dfn}{\stackrel{\triangle}{=}}
\newcommand{\cI}{{\cal I}}
\newcommand{\cJ}{{\cal J}}
\newcommand{\cT}{{\cal T}}
\newcommand{\calF}{{\cal F}}
\newcommand{\calK}{{\cal K}}
\newcommand{\calN}{{\cal N}}
\newcommand{\calX}{{\cal X}}
\newcommand{\EE}{{\Bbb E}}
\newcommand{\complex}{{\Bbb C}}
\newcommand{\ZZ}{{\Bbb Z}}
\newcommand{\lip}{\langle}
\newcommand{\rip}{\rangle}
\newcommand{\one}{\frac{1}{N}\:}
\newcommand{\half}{\frac{1}{2}\:}
\newcommand{\won}{{\mbox{\bf 1}}}
\newcommand{\erfc}{{\rm \,erfc\,}}
\renewcommand{\Im}{{\rm \,Im\,}}
\newcommand{\Prob}{{\rm P\,}}

\newcommand{\Var}{{\rm \,Var\,}}
\newcommand{\bfcdot}{{\mbox{\boldmath$\cdot$}}}
\newcommand{\are}{\mathop{\longrightarrow}_{n\to \infty}}
\newcommand{\areT}{\mathop{\longrightarrow}_{T\to \infty}}

\newcommand{\limsupT}{\limsup_{T \rightarrow \infty}}
\newcommand{\liminfT}{\liminf_{T \rightarrow \infty}}
\newcommand{\limn}{\lim_{n \rightarrow \infty}}
\newcommand{\limsupe}{\limsup_{\ep \rightarrow 0}}
\newcommand{\limsupn}{\limsup_{n \rightarrow \infty}}
\newcommand{\liminfn}{\liminf_{n \rightarrow \infty}}

\newcommand{\ffrac}[2]
  {\left( \frac{#1}{#2} \right)}

\renewcommand{\theequation}{\thesection.\arabic{equation}}

\title{Random polynomials having few or no real zeros}
\date{May 29, 2000}
\author{Amir Dembo
\thanks{Research partially supported by NSF grant DMS-9704552.}
\and
Bjorn Poonen
\thanks{Supported by NSF grant DMS-9801104,
a Sloan Fellowship, and a Packard Fellowship.}
\and
Qi-Man Shao
\thanks{Research partially supported by NSF grant DMS-9802451.}
\and
Ofer Zeitouni
\thanks{Research partially supported by a grant from the Israel Science 
Foundation and by the fund for promotion of research at the 
Technion.
}}

\maketitle

\begin{abstract}
Consider a polynomial of large degree $n$
whose coefficients are independent, identically distributed, nondegenerate
random variables having zero mean and finite moments of all orders. 
We show that such a polynomial has exactly $k$ real zeros
with probability $n^{-b+o(1)}$ as $n \rightarrow \infty$
through integers of the same parity as the fixed integer $k \ge 0$.
In particular, the probability that a random 
polynomial of large even degree $n$ has no real zeros 
is $n^{-b+o(1)}$. The finite, positive constant $b$ 
is characterized via the centered, stationary Gaussian 
process of correlation function 
$\sech (t/2)$. The value of
$b$ depends neither on $k$ nor upon the
specific law of the coefficients. Under an extra smoothness
assumption about the law of the coefficients, with
probability $n^{-b+o(1)}$ one may specify also 
the approximate locations of the $k$ zeros on the real line.
The constant $b$ is
replaced by $b/2$
in case the i.i.d. coefficients have a nonzero mean.
\end{abstract}

\section{Introduction}\label{sec1}
    \setcounter{equation}{0}

Let $\{a_i\}_{i=0}^\infty$ denote a sequence of
independent, identically distributed (i.i.d.)
random variables of zero mean and unit variance.
Consider the random polynomial
\beq{1.0}
f_n(x) = \sum_{i=0}^{n-1} a_i x^i \,.
\eeq
For $n$ odd, define
\beq{1.1}
P_n = \Prob (f_n(x) > 0 \quad \forall x \in \reals)\,.
\eeq
As described in
Section \ref{history},
the study of the number of roots of random polynomials has a long history.
Our main goal is to prove that $P_n=n^{-b+o(1)}$ as $n \rightarrow \infty$
for a finite constant $b>0$,
at least when the coefficient distribution has finite moments of all orders.
The constant $b$ can be described in terms of
the centered stationary Gaussian process $Y_t$ with
correlation function 
$R_y(t) = \sech(t/2)$
(see \req{Ydef} for an explicit construction of $Y_\cdot$).
Define
\beq{1.2}
b = - 4 \lim_{T\to\infty} \frac{1}{T} \log P
\Bigl(\sup_{0\le t\le T } Y_t \le 0\Bigr)  \;,
\eeq
where, throughout this paper, $\log$ denotes the natural logarithm.
The existence of the limit in (\ref{1.2}) and the estimate
$b\in [0.4,2]$ are proved in Lemma \ref{lem-b}. 
We note in passing that our numerical simulations of random polynomials
of degree $n-1 \le 1024$ suggest $b=0.76 \pm 0.03$.

Our main result, which is 
a consequence of Theorem \ref{the-3}
stated in Section \ref{sec-introfull},
is the following
\begin{theorem}
\label{theo-1int}
a) Suppose $\{a_i\}$ is a sequence of zero-mean, unit-variance, i.i.d.
random variables possessing finite moments of all order. Then,
$$
\lim_{n\to\infty}\frac{\log P_{2n+1}}{\log n} = - b\,.$$
b) If $\{a_i\}$ is as above but with $E(a_i)=\mu\neq 0$, we 
denote $P_n^\mu=\Prob (f_n(x) 
\neq
 0 \quad \forall x \in \reals)$.
 Then,
$$
\lim_{n\to\infty}\frac{\log P_{2n+1}^\mu}{\log n} = - b/2\,.$$
\end{theorem}

It is interesting to note that one may answer questions related to 
a prescribed number of zeros. Our main result in this direction is 
the following theorem. For a slightly different variant,
allowing to prescribe the {\it location} of zeros, 
see also Proposition  \ref{prop-2}.
\begin{theo}
\label{theo-2}
Under the assumptions of Theorem \ref{theo-1int} a),
the probability that
the random polynomial $f_{n+1}(x)$ of degree $n$
has $o(\log n/\log \log n)$
real zeros is $n^{-b+o(1)}$ as $n \rightarrow \infty$. 
For any fixed $k$, the probability $p_{n,k}$ that $f_{n+1}$ has exactly $k$
real zeros, all of which are simple, satisfies
$$ \lim_{n\to \infty} 
\frac{\log p_{2n+k,k}}{\log n}
=-b\,.$$
(Obviously, $p_{n,k}=0$ when $n-k$ is odd.)
\end{theo}

The key to our analysis is a detailed study of the case where the
coefficients are Gaussian, implying that 
$f_n(\cdot)$ is a Gaussian 
process (Gaussian processes are particularly 
useful in this context because
for them  comparisons can be made
via Slepian's lemma). The extension to general distribution uses
the strong aproximation results of Koml\'{o}s-Major-Tusn\'{a}dy
\cite{kmt}.
Although this technique requires 
finite moments of all order, we conjecture 
that the asymptotic $n^{-b+o(1)}$ applies to $p_{n,k}$ for $n-k$ even,
whenever the nondegenerate
zero-mean i.i.d. $a_i$ are in the domain of attraction of the
Normal distribution. This conjecture is supported by
the following heuristic derivation of $P_n = n^{-b+o(1)}$.

For $x \in [0,1]$ near $1$, let $x=1-e^{-t}$. Note that
$x^i \approx \exp(-e^{-t} i)$ when $t \gg 0$, and moreover, the function
$h_t(u):=\exp(-e^{-t} u)$ changes slowly in $u$ for $t \gg 0$.
Summation by parts suggests that the sign of $f_n(x)$ is mostly determined by 
the behavior of $\sum_{i=0}^j a_i$ for large $j$ depending on $t$.
Hence,
for $a_i$ in the domain of attraction of the Normal distribution, we
next replace $a_i$ with i.i.d. standard Normal variables $b_i$.
Using the representation $b_i=W_{i+1}-W_i$ for a standard
Brownian motion $W_t$ we further replace the sum over $i=0,\ldots,n-1$
with the corresponding stochastic integral over $[0,\infty)$. This in turn
yields the approximation of the normalized $\hat{f}_n(x):=
f_n(x)/\sqrt{\mbox{\rm Var}(f_n(x))}$ by the centered,
Gaussian process
\beq{Ydef}
Y_t = \frac{\int_0^\infty h_t(u) dW_u}{(\int_0^\infty h_t(u)^2 du)^{1/2}}\;.
\eeq
It is easy to check that the process $Y_\cdot$ of \req{Ydef}
is stationary, with correlation function 
$\sech(t/2)$. By continuity
arguments, $f_n(x)$ typically has a constant sign in $[1-n^{-1},1]$, so
our approximation procedure is relevant only as long as $t \leq \log n$.
Alternatively, $t = \log n$ is where we start seeing $h_t(n)=O(1)$,
contrasting the replacement of the upper limit $n$ in the discrete sum
with the upper limit $\infty$ in the stochastic integral of \req{Ydef}.
We are to consider the possibility of
$f_n(x)=0$ for $x$ in the left and in the right neighborhoods of
both $-1$ and $+1$. In each of these four regimes of $x$ the function
$\hat{f}_n(x)$ is amenable to a similar treatment, leading to an
approximation by the process $Y_t$. With $\hat{f}_n$ having
approximately independent values in the four different regimes, we arrive at
the formula \req{1.2} for $b$.

It is natural to wonder what happens
when $a_i$ are of a symmetric law that is
in the domain of attraction of an $\alpha$-stable law, for some
$\alpha \in (0,2)$. A lower bound on $P_n$ of the form $n^{-c}$
for some finite value of $c$ is then easily obtained by considering
the event
that $a_0$ and $a_{n-1}$ are ``huge'' and positive,
while other coefficients are ``reasonable.''
Repeating the above heuristic for
this case, one is led to believe that the formula \req{1.2}
still applies, but now with $Y_t$ of \req{Ydef} replaced by
\beq{Yadef}
Y_{t,\alpha} = \frac{\int_0^\infty h_t(u) dX^{(\alpha)}_u}
{(\int_0^\infty h_t(u)^\alpha du)^{1/\alpha}}\;,
\eeq
where $X^{(\alpha)}_\cdot$ denotes
the symmetric stable process of index $\alpha$ and the stochastic
integral in \req{Yadef} is to be interpreted via integration by parts.
We have yet no strong evidence to support the above statement. However,
our numerical simulations indicate the behavior
$P_n = n^{-b_1 + o(1)}$ for i.i.d. Cauchy random variables $a_i$
(that is, $\alpha=1$), where $b_1 \approx 0.86$ is larger than $b$.

\subsection{Historical remarks} \label{history}

Throughout this section, $\{a_i\}$ are
independent, identically distributed,
nondegenerate, real-valued random variables.\footnote{
    Some authors whose work we mention assumed $a_0=1$ or $a_n=1$,
    but as far as asymptotic behavior as $n \rightarrow \infty$
    is concerned, it makes little difference.}
Let $N_n$ denote the number of distinct\footnote{
    The asymptotic behavior does not depend on whether
    roots are counted with multiplicity or not.}
real zeros of $f(x):=\sum_{i=0}^n a_i x^i$.
(For the sake of definiteness,
we define $N_n=0$ when $f$ is the zero polynomial.)
So, $p_n:=P(N_n=0)=p_{n,0}$ and we also
let $E_n$ and $V_n$ denote the mean and variance of $N_n$.

The study of real roots of random polynomials has a long and full history,
but most of it deals with the asymptotic behavior of $E_n$
instead of $p_n$.
Presumably this is because $E_n$ is
much easier to estimate: because expectation is linear,
one can compute $E_n$ by integrating over the real line the probability
of having a root in $(t,t+dt)$, for example.

Although as mentioned in~\cite[p.~618]{todhunter},
one can find probabilistic statements in the context of roots of polynomials
as early as 1782 (Waring) and 1864 (Sylvester),
the first people to study the asymptotic behavior of $N_n$ seem to be
Bloch and P\'olya~\cite{blochpolya}.
In 1932, they proved $E_n=O(n^{1/2})$ for the coefficient distribution
$P(a_i=1)=P(a_i=0)=P(a_i=-1)=1/3$.
This work led Littlewood and Offord to undertake a systematic study of $N_n$
in a series of papers \cite{lo1},\cite{lo2},\cite{lo3}
starting in 1938.
They proved that if
the $a_i$ are all uniform on $[-1,1]$,
or all Normal, or all uniform on $\{-1,1\}$, 
(i.e. $P(a_i=1)=P(a_i=-1)=1/2$),
then
$$
P \left(N_n>25(\log n)^2 \right) \le \frac{12 \log n}{n},
\quad \mbox{and} \quad
P\left( N_n< \frac{\alpha \log n}{(\log \log n)^2}  \right)
< \frac{A}{\log n}
$$
for some constants $\alpha$ and $A$.
In particular, for some constant $\alpha'$,
$$
\frac{\alpha' \log n}{(\log \log n)^2}
\le E_n \le 25(\log n)^2 + 12 \log n
$$
and $p_n=O(1/\log n)$ for these distributions.
This upper bound for $p_n$ has apparently not been improved,
until the current paper.\footnote{
The only result in the literature that might be said
to have improved our knowledge of $p_n$
is~(\ref{variance}),
which together with~(\ref{enfornormalattraction})
implies for many distributions
that $\limsup_{n \rightarrow \infty} p_n \log n \le \pi - 2$.
The bound has the same form as that arising from the work
of Littlewood and Offord,
but the constant has been made explicit.}

In 1943 Kac~\cite{kac1} found the exact formula
\beq{kacexact}
E_n ={ 1 \over \pi}
\int_{-\infty}^\infty
\sqrt{ { 1 \over (t^2-1)^2}
- { (n+1)^2 t^{2n} \over (t^{2n+2}-1)^2} } dt  \;,
\eeq
when $a_i$ is Normal with mean zero,
and extracted from it the asymptotic estimate
\begin{equation}
\label{enfornormalattraction}
    E_n \sim \frac{2}{\pi} \log n \;.
\end{equation}
Much later Jamrom~\cite{jamrom} and Wang~\cite{wang}
improved this to $E_n = (2/\pi) \log n + C + o(1)$
for an explicit constant $C$,
and ultimately Wilkins~\cite{wilkins2} obtained
an asymptotic series for $E_n$ from \req{kacexact}.
In 1949 Kac~\cite{kac2} obtained \req{enfornormalattraction}
for the case where $a_i$ is uniform on $[-1,1]$.
Erd\H{o}s and Offord~\cite{erdosofford}
obtained the same asymptotics for $a_i$ uniform on $\{-1,1\}$.
Stevens~\cite{stevens} proved \req{enfornormalattraction}
for a wide class of distributions, and this estimate was
finally extended by Ibragimov and Maslova~\cite{im1},\cite{im2}
to all mean-zero distributions in the domain of attraction
of the Normal law.

At around the same time (the late 1960's),
Logan and Shepp~\cite{loganshepp1},\cite{loganshepp2}
discovered that if the coefficient distribution is 
the symmetric stable distribution with characteristic function
$\exp(-|z|^\alpha)$, $0 < \alpha \le 2$, then $E_n \sim c_\alpha \log n$,
where
    $$c_\alpha := \frac{4}{\pi^2 \alpha^2}
    \int_{-\infty}^{\infty} dx\, \log \int_0^\infty
    \frac{|x-y|^\alpha e^{-y}}{|x-1|^\alpha} \, dy
> \frac{2}{\pi}.$$
They also proved $\lim_{\alpha \rightarrow 0^+} c_\alpha = 1$,
and performed calculations that suggested
that $c_\alpha$ is a decreasing function of $\alpha$,
terminating at $c_2=2/\pi$, Kac's value for the Normal distribution.
Ibragimov and Maslova~\cite{im4}
extended these results by finding the asymptotic behavior of $E_n$
for arbitrary distributions in the domain of attraction
of a stable distribution.
The asymptotics are different when the distribution has nonzero mean;
for instance~\cite{im3},
if $a_i$ are Normal with nonzero mean,
then $E_n \sim (1/\pi) \log n$ instead of $(2/\pi) \log n$.
Shepp (private communication) has conjectured that there exists
a universal constant $B$ such that
    $$\limsup_{n \rightarrow \infty} \frac{E_n}{\log n} \le B$$
for {\em any} coefficient distribution
(satisfying only the hypotheses at the beginning of this section).
If $B$ exists, then $B \ge 1$ by the work of Logan and Shepp
mentioned above.

In 1974, Maslova~\cite{maslovavariance},\cite{maslovanormal}
proved that if $P(a_i=0)=0$, $Ea_i=0$
and $E(a_i^{2+\epsilon}) < \infty$ for some $\epsilon>0$,
then
\begin{equation}
\label{variance}
    V_n \sim \frac{4}{\pi} \left(1 - \frac{2}{\pi} \right) \log n
\end{equation}
and $N_n$ is asymptotically Normal.

Much work was also done on complex roots of $f_n$; see \cite{IZ}
and references therein for an updated account.
Further results on random polynomials and their generalizations
can be found in the books~\cite{bharuchareid,fara}
and the survey article~\cite{edelmankostlan}.

Our interest in the asymptotics of $p_n$
grew out of a problem in arithmetic geometry.
The paper~\cite{poonenstoll} showed that
Jacobians of curves over $\Q$ could be odd,
in the sense of having Shafarevich-Tate groups of nonsquare order
(despite prior claims in the literature that this was impossible).
Moreover it was shown (in a sense that was made precise)
that the probability
that a random hyperelliptic curve $y^2=f(x)$ of genus $g$ over $\Q$
has odd Jacobian
could be related to a sequence of ``local'' probabilities,
one for each nontrivial absolute value on $\Q$.
The computation of the local probability for the standard
archimedean absolute value reduced to the knowledge of
the probability that the curve $y^2=f(x)$ has no real point,
or equivalently,
the probability that the random polynomial $f(x)$
satisfies $f(x)<0$ for all real $x$.
Although the asymptotic behavior of this probability
was not needed in a substantial way in~\cite{poonenstoll},
the authors of that paper found the question to be of sufficient
interest in its own right
that they developed heuristics that led them to conjecture
the existence of a universal constant $b>0$
such that $p_n=n^{-b+o(1)}$,
for any mean-zero distribution
in the domain of attraction of the Normal law.

\subsection{Statement of main theorems}
\label{sec-introfull}

Let $\hat{f}_n(x):=f_n(x)/\sqrt{E(f_n(x)^2)}$ denote the normalized
random polynomial, so $\hat{f}_n(x)$ has unit variance for each $x$.
Instead of proving only $P_n=n^{-b+o(1)}$,
we generalize in the following, to facilitate applications to related problems.
\begin{theo}
\label{the-3}
Suppose $a_i$ are zero-mean i.i.d.
random variables
of unit variance and
with finite moments of all orders.
For $n-1$ even, let
$$
P_{n,\gamma_n} = P \Bigl(\hat{f}_n(x) > \gamma_n(x)
\quad \forall x \in \reals \Bigr) \;,
$$
for nonrandom functions $\gamma_n(x)$ such that
$n^{\delta}|\gamma_n(x)| \to 0$ uniformly in $x \in \reals$ for some
$\delta>0$. Then,
\beq{1.3}
\frac{\log P_{n,\gamma_n}}{\log n} \are - b \;.
\eeq
The upper bound $P_{n,\gamma_n} \leq n^{-b+o(1)}$ applies as soon as
$$
\inf\{ \gamma_n(x) : ||x|-1|\leq n^{-\ep_n} \} \to 0 \quad\quad
{\rm for \ any \ } \ep_n \to 0 \;.
$$
\end{theo}

The key to the proof of Theorem \ref{the-3} is the analysis of
$P_{n,\gamma_n}$ for random polynomials $f_n(x)$
with coefficients $\{b_i\}$ that are i.i.d. standard Normal variables.
To distinguish this case, we use 
throughout the notations $f^b_n(x)$, $\hat{f}^b_n$
and $P^b_{n,\gamma_n}$ for $f_n(x)$, $\hat{f}_n(x)$ and $P_{n,\gamma_n}$,
respectively, when dealing with polynomials of coefficients that are
Normal variables. The next theorem summarizes our results in this special case.
\begin{theo}
\label{the-1}
The convergence of $\log P^b_{n,\gamma_n}/\log n$ to $-b$ applies in the
standard Normal case, as soon as
the nonrandom functions $\gamma_n(x) \leq M < \infty$ are such that
$$
\sup\{ |\gamma_n(x)| : ||x|-1|\leq n^{-\ep_n} \} \to 0 \quad\quad
{\rm for \ some \ } \ep_n \to 0 \;.
$$
\end{theo}

The following proposition is the variant of
Theorem \ref{theo-2} alluded to above. 
It shows that
with
probability $n^{-b +o(1)}$
one may also prescribe arbitrarily the location of the $k$
real zeros of $f_{n+1}(x)$, provided the support of the law of $a_i$
contains an open interval around $0$.
The latter assumption is to some extent necessary. For example, when
$P(a_i=1)=P(a_i=-1)=1/2$ it is easy to see that
$f_{n+1}(x)$ cannot have zeros in $[-1/2,1/2]$.
\begin{prop}
\label{prop-2}
Suppose $a_i$ are zero-mean i.i.d.
random variables
of unit variance,
finite moments of all orders, and the support of the law of each $a_i$
contains the interval $(-\eta,\eta)$ for some $\eta>0$.
Given disjoint open intervals $U_1,\ldots,U_\ell$ and positive
integers $m_1,\ldots,m_\ell$, the probability that
the random polynomial $f_{n+1}(x)$ has exactly $m_i$ real zeros in $U_i$
for each $i$ and no real zeros anywhere else is $n^{-b+o(1)}$
for $n \to \infty$ through integers of the same parity as $k=\sum_i m_i$.
\end{prop}

The organization of this paper is as
follows.
Auxiliary lemmas 
about Gaussian processes,
needed for the proof of Theorem~\ref{the-1},
are grouped in Section~\ref{sec2} (including the bounds on $b$ mentioned in
the introduction, c.f.
Lemma \ref{lem-b}).
Relying upon Gaussian techniques,
the proof of the lower bound of Theorem~\ref{the-1}
is in Section~\ref{sec3}, and the complementary upper bound in
Section~\ref{sec4}. Building upon Theorem \ref{the-1},
and with the help of strong approximation,
Section \ref{sec7} provides the proof of our main result, Theorem \ref{the-3}.
Theorem \ref{theo-1int}
is
then derived in Section \ref{sec4new}.
Section \ref{sec6} provides the upper bound on the
probability of interest in Theorem \ref{theo-2}, with
the lower bound proved in Section \ref{sec5}.
Finally, Proposition \ref{prop-2} is proved in Section \ref{sec10}.

\section{Auxiliary lemmas} \label{sec2}

            \setcounter{equation}{0}

We start by introducing several notations that appear throughout this
work. For $n$ odd, let $c_n(x,y)$ denote the covariance function of
$\hat{f}_n(x)$, that is
\beq{1.4}
c_n(x,y) = \frac{E\Bigl(f_n(x) f_n(y) \Bigr)}
{\sqrt{E(f_n(x)^2) E(f_n(y)^2)}}
\eeq
Then, for $x \neq \pm 1$ and $y \neq \pm 1$,
\begin{subequations}
\label{1.5}
\beq{1.5a}
c_n(x,y) = \frac{g(x^n, y^n)}{g(x,y)}
\eeq
where
\beq{1.5b}
g(x,y) = \frac{|xy-1|}{\sqrt{|(1-x^2)(1-y^2)|}} \ge 0 \,.
\eeq
Note that
$g(x,y) = g(-x,-y) = g(\frac{1}{x},\frac{1}{y})$.
Further,
\beq{1.5c}
\forall x, y\in
(-1, 1),  \quad g(x,y) \ge 1
\eeq
\end{subequations}
and the change of variables $z=1-x, w=1-y$, leads to
\beq{1.6}
\frac{1}{g(x,y)}
= \frac{2\sqrt{zw}}{z+w}
\Bigl[1-\Bigl[\frac{1-\frac{wz}{z+w} - \sqrt{1-\frac{z}{2}}\;
\sqrt{1-\frac{w}{2}}}{1-\frac{wz}{z+w}} \Bigr]\Bigr]
\eeq

A good control on $g(x,y)$ is provided by the following lemma.
\begin{lemma}
\label{lem-0}
For any $z, w \in (0,1/2]$
$$
\frac{1}{8} (w-z)^2 \le
\Bigl(1-\frac{wz}{w+z} - \sqrt{1-\frac{z}{2}}\;
\sqrt{1-\frac{w}{2}}\Bigr) \frac{\max(z,w)}{\Bigl(1-\frac{wz}{w+z}\Bigr)}
\le
(w-z)^2 \,.
$$
\end{lemma}

\proof   Let $z+w= \eta$, $z-w=\xi$, assuming without loss of
generality that $0 < \xi\le \eta \le 1$. Since
$$
1\ge 1 - \frac{wz}{w+z}
\ge \half\,,\quad
1 \ge \frac{\max (z,w)}{(z+w)} \ge \half\,,
$$
it suffices to prove that
$$
f(w,z) = \frac{1}{(z-w)^2} \Bigl[
1-\frac{wz}{w+z} - \sqrt{1-\frac{z}{2}}\;\sqrt{1-\frac{w}{2}}\;\Bigr]
(z+w) \in [1/4, 1/2] 
 \;.
$$
To this end, observe that for all $0 < \xi \le \eta \le 1 $
we have
\beaa
f(w,z) & = &
\frac{\eta}{\xi^2} \Bigl[ 1 - \frac{\eta^2 - \xi^2}{4\eta} -
\sqrt{\Bigl(1-\frac{\eta}{4}\Bigr)^2 - \Bigl(\frac{\xi}{4}\Bigr)^2}
\:\Bigr] \\
& = & \frac{\eta}{\xi^2} \Bigl[\Bigl(1-\frac{\eta}{4} \Bigr) +
\frac{\xi^2}{4\eta} - \sqrt{\Bigl(1-\frac{\eta}{4}\Bigr)^2 -
\Bigl(\frac{\xi}{4}\Bigr)^2}\: \Bigr]
 = \frac{1}{4} + \frac{\tilde{\eta}}{4 \tilde{\xi}^2} \Bigl[ 1-
\sqrt{1-\tilde{\xi}^2} \Bigr] \,,
\eeaa
where $\tilde{\xi} = \xi/(4-\eta)$, $\tilde{\eta} = \eta/(4-\eta)$.
Since $\xi\le \eta \le 4 -\eta$
 and $0 \le 1 - \sqrt{1-\tilde{\xi}^2} \le \tilde{\xi}^2$ it
follows that
$1/4 \le f(w,z) \le 1/4 + \tilde\eta/4 \le 1/2$ as needed.
\qed

The control of Lemma \ref{lem-0} on $g(x,y)$, hence on $c_n(x,y)$,
shall give rise to the perturbed centered Gaussian processes
$Y^{(\alpha)}$ of the next lemma.
\begin{lemma}
\label{lem-2.2}
Let $\alpha \in [0,1]$ and define the covariance
$$
R^{(\alpha)}(\tau) = 
\sech \Bigl({\tau}/{2}\Bigr)
\Bigl\{1-\alpha (1-e^{-|\tau|})^2 \Bigr\} \,.
$$
Then there exist independent, stationary centered
Gaussian processes $Y_t, Z_t$, with covariances $R_y(\tau)=R^{(0)}(\tau)$ and
$$
R_z(\tau) = R^{(1)}(\tau) = 
\sech \Bigl({\tau}/{2}\Bigr)
\Bigl(2e^{-|\tau|} - e^{-2|\tau|} \Bigr)
$$
respectively, such that the process
$Y_t^{(\alpha)} := \sqrt{1-\alpha}\: Y_t + \sqrt{\alpha}Z_t$
has covariance $R^{(\alpha)}(\tau)$.
\end{lemma}

\proof
Since $R^{(\alpha)}(\tau) = (1-\alpha) R_y (\tau) + \alpha R_z(\tau)$, all one
needs is to check that both $R_y(\tau)$ and $R_z(\tau)$ are
covariance functions,
i.e. to check that their Fourier transforms are nonnegative.  To this end,
note that
$$
S_y(\om):= \calF (R_y(\tau)) =
\int_{-\infty}^\infty e^{i \om \tau} R_y(\tau) d\tau \\
 = 2 \int_0^\infty 
\cos (\om\tau)\sech(\tau/2) d\tau = 2\pi\sech (\om\pi) \ge 0\,,
$$
c.f. \cite[p.~503, formula 3.981.3]{GR}.
Furthermore,
$$
S_z(\om):=\calF (R_z(\tau)) =
\int_{- \infty}^{\infty} e^{i\om \tau} R_z (\tau) d\tau = S_y (\om) \ast
F(\om)
\,,
$$
where $\ast$ stands throughout for the convolution operation and
\beaa
F(\om) & = & \int_{-\infty}^\infty e^{i\om\tau} \Bigl( 2e^{-|\tau|} -
e^{-2|\tau|}\Bigr) d\tau
 =  2 \int_0^\infty \cos (\om\tau) \Bigl(2e^{-\tau} -e^{-2\tau}\Bigr)
d\tau \\
& = & \frac{4}{1+\om^2} - \frac{4}{4+\om^2} = \frac{12}{
(4+\om^2)(1+\om^2)}
\ge 0
\,.
\eeaa
Hence, $S_z  (\om) \ge 0$.
\qed

The effect of  nonrandom functions $\gamma_n(x)$ as well as that of
considering the processes $Y^{(\alpha)}$ for some $\alpha_n \downarrow 0$
are dealt with by the continuity properties of $Y_t$ and $Z_t$ outlined
in the next lemma.
\begin{lemma}
\label{lem-2.3}
Let $Y_t,Z_t$ be as in Lemma~\ref{lem-2.2}. Then, for any
positive $\ep_T \to 0$,
\beq{2.5}
\lim_{T \to \infty}
P \Bigl(\sup_{0\le t \le T} Z_t \leq \sqrt{\ep_T^{-1} \log T} \Bigr) = 1\,,
\eeq
whereas
\beqn{2.6}
\limsupT \frac{1}{T} \log P \Bigl(\sup_{0\le t \le T} Y_t \le \ep_T\Bigr)
&=&
\liminfT \frac{1}{T} \log P \Bigl(\sup_{0\le t \le T} Y_t \le
-\ep_T\Bigr)
\nonumber \\
&=&
\lim_{T\to\infty} \frac{1}{T} \log P \Bigl(\sup_{0\le t \le T} Y_t \le 0
\Bigr)
= - \frac{b}{4} \,.
\eeqn
Moreover, for any positive $\gamma_T \to 0$ and $\alpha_T \log T \to 0$,
\beq{2.6n}
\liminf_{T\to\infty} \frac{1}{T} \log P \Bigl(\inf_{0\le t \le T}
Y_t^{(\alpha_T)} \ge \gamma_T \Bigr) \geq - \frac{b}{4} \,.
\eeq
\end{lemma}

\proof 
The existence of the limit in the right hand side of \req{2.6}
(and hence in \req{1.2}) is ensured by sub-additivity: since
$R_y(\cdot) > 0$, Slepian's lemma (c.f. \cite[Page 49]{adler}),
and the stationarity of $Y_\cdot$ imply
\beaa
P\Bigl( \sup_{0\le t\le T+S} Y_t\le 0 \Bigr) & \ge &
P\Bigl(\sup_{0\le t \le T} Y_t \le 0\Bigr)
P\Bigl(\sup_{T\le t \le T+S } Y_t \le 0\Bigr) \\
& = &
P\Bigl(\sup_{0\le t \le T} Y_t \le 0\Bigr)
P\Bigl(\sup_{0\le t \le S} Y_t \le 0\Bigr) \,.
\eeaa

Fix $\ep_T \to 0$. From Lemma~\ref{lem-2.2}, we have that
$$
S_z(\om) = \calF(R_z (\tau)) =
2\pi \sech (\om \pi) \ast
\frac{12}{(4+\om^2) (1+\om^2)}\,,
$$
which implies that $\sup_\om \{S_z(\om) \om^4\} < \infty$.
Hence $\int_{-\infty}^{\infty} \om^2 S_z(\om)d\om < \infty$.
It follows that
$$
-\frac{\partial^2}{\partial \tau^2}
R_z(\tau)\Big|_{\tau=0} = E(\dot{Z}_t^2) < \infty \,.
$$
Since $|Z_t| \le |Z_0| + \int_0^1 |\dot{Z}_t| dt$, it follows by
stationarity of the centered Gaussian process $\dot{Z}_t$ that
$$
m_1 := E (\sup_{0 \le t \le 1} |Z_t|) \le
\sqrt{E (Z_0^{2})} + \sqrt{E(\dot{Z}_t^2)} < \infty \,.
$$
By the stationarity of $Z_t$
and Borell's inequality (c.f. \cite[Page 43]{adler}),
for all $\lam \geq m_1$,
$$
P \Bigl( \sup_{0\le t \le T} Z_t \ge \lam \Bigr)
 \le T P\Bigl(\sup_{0\le t\le 1} |Z_t|\ge \lam \Bigr)
\le 2 T  \exp\Bigl(-\frac{(\lam-m_1)^2}{2 R_z(0)}\Bigr) \,.
$$
Setting
$\lambda=\sqrt{\ep_T^{-1} \log T}$
we obtain that as $T \to \infty$,
$$
P \Bigl( \sup_{0\le t \le T} Z_t \ge \sqrt{\ep_T^{-1} \log T} \Bigr) \to 0 \;,
$$
which yields \req{2.5}.

To see \req{2.6}, let the Gaussian law of $Y_\cdot$
on $C(\reals;\reals)$ be denoted
by $P_y$. Let ${\cal R}_y$ denote the covariance operator
associated with $P_y$, that is,
${\cal R}_y g(t) = \int_{-\infty}^{\infty} {\cal R}_y (t-s) g(s) ds$, with
$K_y = {\cal R}_y^{-1}$ denoting its inverse (defined on the range of
${\cal R}_y$). We also let $\lip \cdot, \cdot \rip$ denote the
inner product of $L^2(\reals)$. Fixing $T<\infty$ note that the
deterministic function
$f_T(t):= \ep_T \exp(\frac{1}{4}-\Bigl(\frac{t}{T}\Bigr)^2)$ is
in the Reproducing Kernel Hilbert Space (RKHS) associated with the process
$Y_\cdot$. Indeed, the Fourier transform of $f_T$
is $\hat f(\om) = c_1  T \ep_T e^{-\om^2 T^2}$ (for some
$c_1<\infty$), so it follows by Parseval's theorem that for some $c_2<\infty$
and all $T$,
\beq{2.9}
\lip f_T, K_y f_T\rip  =
\int_{-\infty}^{\infty} \frac{| \hat f (\om)|^2}{S_y(\om)} d\om \\
 = 2\pi \int_{-\infty}^{\infty} 
|\hat f (\om)|^2 \sech(\om \pi)
d\omega \leq c_2 T \ep_T^2 \,,
\eeq
In particular, $\lip f_T, K_y f_T\rip$ is finite and
the Radon-Nikod\'ym derivative
$$
\Lambda_T(Y) = \exp \Bigl(\lip f_T, K_y Y\rip - \half
\lip f_T, K_yf_T \rip \Bigr) \;,
$$
is well defined and finite for $P_y$-almost-every $Y_\cdot$.
Since
$f_T(t) \ge \ep_T$ for all $-\frac{T}{2} \le t \le \frac{T}{2}$,
it follows that
\beqn{2.7}
&&
P \Bigl( \sup_{- \frac{T}{2} \le t \le \frac{T}{2}} Y_t \le \ep_T\Bigr)
= P\Bigl(\sup_{-\frac{T}{2}\le t \le \frac{T}{2}} \{ Y_t-\ep_T \}
\le 0 \Bigr)
 \leq P\Bigl(\sup_{-\frac{T}{2}\le t \le \frac{T}{2}} \{ Y_t- f_T(t) \}
\le 0 \Bigr)
\nonumber \\
& = & E \Bigl( \Lambda_T (Y) \,
\won_{ \{\sup_{-\frac{T}{2}\le t \le \frac{T}{2}} Y_t
\le 0 \} } \Bigr)
\le E (\La_T(Y)^q)^{\frac{1}{q}} \Bigl[P
\Bigl(\sup_{- \frac{T}{2} \le t \le
\frac{T}{2}} Y_t \le 0 \Bigr)\Bigr]^{\frac{1}{p}}
\eeqn
where $\frac{1}{p} + \frac{1}{q} = 1$. Note that
$$
\Bigl(E(\La_T(Y)^q)\Bigr)^\frac{1}{q} = \exp \Bigl(\frac{q-1}{2}
\lip f_T, K_y f_T \rip \Bigr) \,.
$$
Hence, choosing $q_T =(1/\ep_T) \to \infty$ it follows from \req{2.9} that
$$
\frac{1}{T} \log \Bigl( E (\La_T(Y)^q) \Bigr)^{\frac{1}{q}} \areT 0 \,.
$$
Substituting in \req{2.7} and using the stationarity of $Y_\cdot$ and
existence of the limit 
in \req{1.2}, one has that
\beq{2.8}
\limsupT \frac{1}{T} \log P \Bigl(\sup_{0 \le t \le T} Y_t \le \ep_T
\Bigr) \le \lim_{T\to \infty} \frac{1}{T} \log
P\Bigl(\sup_{0 \le t \le T} Y_t \le 0 \Bigr) \,.
\eeq
The equality in \req{2.8} is
then obvious. The other equality in \req{2.6} follows by a similar proof,
starting with
$$
P\Big( \sup_{- \frac{T}{2} \le t \le \frac{T}{2}} Y_t \le 0\Big)
\leq
P\Big(\sup_{-\frac{T}{2}\le t \le \frac{T}{2}} \{ Y_t- f_T(t) \}
\le - \ep_T \Big)
\leq
E (\La_T(Y)^q)^{\frac{1}{q}}
P\Big( \sup_{- \frac{T}{2} \le t \le
\frac{T}{2}} Y_t \le - \ep_T \Big)^{\frac{1}{p}}.
$$

Turning to prove \req{2.6n},
set $\ep_T = 3 \max( \gamma_T, (\al_T \log T)^{1/3}) \to 0$, and note that
$$
\sqrt{1-\al_T} \ep_T - \sqrt{\al_T}
 \sqrt{\ep^{-1}_T \log T }\geq \gamma_T \;,
$$
once $T$ is large enough that $\al_T \leq 1/3$.
Then, by the independence of $Y_t$ and $Z_t$,
$$
P \Bigl(\inf_{0\le t \le T} Y_t^{(\alpha_T)} \ge \gamma_T \Bigr) \geq
P \Bigl(\inf_{0\le t \le T} Y_t \ge \ep_T \Bigr)
P \Bigl(\inf_{0\le t \le T} Z_t \ge -\sqrt{\ep_T^{-1} \log T} \Bigr) \;.
$$
With the laws of the processes $Y_t$ and $Z_t$ invariant to a change of sign,
the inequality \req{2.6n} is thus a direct consequence of \req{2.5} and
\req{2.6}.
\qed

The control of $\hat{f}^b_n(x)$ for $x \in [1-n^{-1},1]$ is achieved in the
next lemma by means of the sample path smoothness of $f^b_n(\cdot)$.
\begin{lemma}
\label{lem-2.1}
For any finite $\gamma$, the set of limit points of
$C_n=P(\hat{f}^b_n(x) > \gamma, \, \forall x\in [1-n^{-1}, 1]\,)$
is bounded below by some $C_\infty = C_\infty(\gamma) >0$.
\end{lemma}

\proof Without loss of generality we assume that $\gamma \geq 0$.
Since $x \mapsto E(f^b_n(x)^2)$ is increasing on $[0,\infty)$, with
$E(f^b_n(1)^2)=n$, it follows that for any $\la> 0$
\beqn{2.1}
C_n &\ge&
P\Bigl(f^b_n(x) > \gamma \sqrt{n}, \quad \forall x\in [1-n^{-1}, 1]\, \Bigr)
\nonumber \\
&\ge& P \Bigl(f^b_n (1) > (\la+\gamma) \sqrt{n} \Bigr)
- P\Bigl(
\sup_{(1-n^{-1}) \le \xi \le 1} \fbnp (\xi) \ge \la n^{3/2}\Bigr)
\eeqn
We wish to apply Borell's
inequality to bound the second term in \req{2.1}.  To this end, note
that
$$
\fbnp (\xi) = \sum_{i=0}^n i b_i \xi^{i-1}
 =  \sum_{i=0}^n \Bigl(\sum_{j=0}^i b_j \Bigr)
\Bigl[ i \xi^{i-1} - (i+1) \xi^i\Bigr] + (n+1) \xi^n
\sum_{j=0}^n b_j \,.
$$
By Kolmogorov's maximal inequality,
$$
E \Bigl[\sup_{i\le n} \Bigl| \sum_{j=0}^i b_j \Bigr| \Bigr]
\le c_1 n^{1/2} \,.
$$
Hence, for some $c_2>0$,
\beq{2.3}
E \Bigl| \sup_{1-n^{-1} \le \xi \le 1} \fbnp (\xi) \Bigr|
\le c_2 n^{{3}/{2}} \,.
\eeq
Furthermore, we have that
$$
\sup_{1-n^{-1} \le \xi \le 1}
E \Bigl[\frac{1}{n^{3/2}} \fbnp (\xi) \Bigr]^2
= \frac{1}{n^3} \sum_{i=1}^{n-1} i^2 \are \frac{1}{3} \,,
$$
implying, by Borell's inequality and \req{2.3}, that for some
finite $c_3$, all $n$ and any $\la \geq c_2$,
\beq{2.4}
P\Bigl(\sup_{1-n^{-1} \le \xi \le 1} \fbnp  (\xi) \ge \la n^{3/2} \Bigr)
\le c_3 e^{-3 (\la-c_2)^2/2} \;.
\eeq
Since $n^{-1/2} f^b_n (1)$ is a standard Normal random variable,
it follows that for some positive $c_4=c_4(\gamma)$, $\la=\la(\gamma)$
large enough and all $n$,
\beq{2.2}
P\Bigl(f^b_n (1) > (\la+\gamma) \sqrt{n} \Bigr) \ge c_4 e^{-\la^2}
\ge 2 c_3 e^{-3 (\la-c_2)^2/2} \;.
\eeq
Substituting \req{2.4} and \req{2.2} in \req{2.1}, one concludes that
$\liminfn C_n \geq C_\infty >0$ as claimed.
\qed

The next lemma provides the bounds on the value of $b$ stated in
the introduction.
\begin{lemma}
\label{lem-b}
The limit in \req{1.2} exists, and the
constant $b$ there satisfies the bounds $0.4 \leq b \leq 2$.
\end{lemma}

\proof
The existence of the limit in \req{1.2}
was proved 
in the course of proving Lemma \ref{lem-2.3}.
Recall that $R_y(t) \geq e^{-|t|/2}$, the covariance of the
stationary Ornstein-Uhlenbeck process $X_\cdot$
As can be checked by computing the covariance,
a  representation of the process $\{X_t\}$ can be obtained as
\beq{o-u}
X_t=e^{-t/2}V_{e^t}=e^{-t/2}(V_{e^t}-V_1+V_1)=
e^{-t/2}(W_{e^t-1}+X_0)\,,
\ee
for some standard Brownian motions 
$V_\cdot$, $W_\cdot$ and a standard normal random variable $X_0$
that is independent of $W_\cdot$. Hence, for $\eta=(e^T-1)^{-1/2}$,
\begin{eqnarray*}
P(\sup_{0\le t\le T } X_t \le 0)&=&
E[1_{X_0 \leq 0}\, P(\sup_{0\leq t\le e^T-1} \{W_t\} \leq -X_0 |X_0)] \\
&=& E[1_{X_0 \leq 0}\, (1-2P(W_{e^T-1} \geq -X_0|X_0))] \\
&=& \frac{1}{\pi} \int_{-\infty}^0 \int_{0}^{-\eta x} e^{-(x^2+y^2)/2} dy dx
= \frac{1}{\pi} \arctan(\eta) = \frac{1}{\pi} e^{-T/2} (1+o(1)).
\end{eqnarray*}
Consequently, Slepian's lemma implies the bound $b \leq 2$. 

The proof of the complementary bound is based on the following observation.
Suppose that ${\bf X} \in \reals^n$ and ${\bf Y}\in\reals^n$ are zero-mean,
normally distributed random vectors 
with covariance matrices $\Sigma_x$ and $\Sigma_y$ respectively.
If $\Sigma_x - \Sigma_y$ is a positive semidefinite matrix, then
the Radon-Nikod\'ym derivative of the law of ${\bf Y}$ with respect to that
of ${\bf X}$ is at most
$( \det \Sigma_x / \det \Sigma_y \Big)^{1/2}$, hence
\beq{b2}
P({\bf Y} \in C)  \leq
\Big( { \det \Sigma_x \over \det \Sigma_y } \Big)^{1/2} P({\bf X} \in C),
\eeq
for all $C \subset \reals^n$ (c.f. \cite[Lemma 3.1]{shao}).
\ignore{
{\it Proof.} Let
$f_X$  and
$f_Y$  be the joint density functions of $\X$ and $\Y$,
respectively.
Since $\Sigma_X-\Sigma_Y$ is positive semidefinite,
$\Sigma_Y^{-1} - \Sigma_X^{-1} $ is positive semidefinite, too.
Hence
\beaa
f_Y(\x)
& = & { 1 \over (2\pi)^{n/2} |\Sigma_Y|^{1/2}}
\exp\Big(- { 1\over 2} \x' \Sigma_Y^{-1} \x \Big)\\
&\leq &  { 1 \over (2\pi)^{n/2} |\Sigma_Y|^{1/2}}
\exp\Big(- { 1\over 2} \x' \Sigma_X^{-1} \x \Big)
= (|\Sigma_X| /|\Sigma_Y|)^{1/2} f_X (\x),
\eeaa
which yields (\ref{b2}) immediately. \qed
}
Indeed, to prove that $b \geq 0.4$, it suffices to show that
\beq{b3}
P(\max_{1 \leq i \leq n} Y_{5i} \leq 0)
\leq \exp(- 0.5\, n)
\eeq
for all $n \geq 2$. Let 
$$\rho= 2\, \sech (5/2)= {2 e^{-2.5} \over 1+e^{-5}},  \
 \ \la_0= 4\sum_{i=2}^\infty \sech(5i/2)=
4 \sum_{i=2}^\infty {e^{-2.5i} \over 1+ e^{-5i}} \
$$
and $(X_1, \cdots, X_n)$ be independent normal random variables
each having zero mean and variance $\la:=1+2\rho+\la_0$.
Denote the covariance matrices
of $(X_i, 1 \leq i \leq n)$ and $(Y_{5i}, 1 \leq i \leq n)$
by $\Sigma_x$ and $\Sigma_y$, respectively.
It is easy to see that $\Sigma_x- \Sigma_y$ is a dominant
principal diagonal matrix and as such it is positive semidefinite.
Thus, by (\ref{b2})
$$
P(\max_{1 \leq i \leq n} Y_{5i} \leq 0)
\leq \Big( {\det \Sigma_x \over \det \Sigma_y}\Big)^{1/2}
P(\max_{1 \leq i \leq n} X_i \leq 0) 
={ \la^{n/2} 2^{-n}\over (\det \Sigma_y)^{1/2} }.
$$
To estimate $\det \Sigma_y$, let
$\Sigma_n=(r_{ij}, 1 \leq i, j \leq n)$ be a tri-diagonal
matrix with $r_{ii}= 1-\la_0$, $r_{i,i+1}=r_{i-1,i}=\rho$
and $r_{ij}=0$ for other $i,j$.
Then, $\Sigma_y - \Sigma_n$ is a positive semidefinite matrix and hence
$$
\det \Sigma_y \geq \det \Sigma_n := D_n.
$$
Since $D_n=(1-\la_0) D_{n-1}-\rho^2 D_{n-2}$, direct calculation shows that
$$
D_n \geq \Big( {1 \over 2} (1-\la_0 + \sqrt{(1-\la_0)^2 - 4 \rho^2})\Big)^n .
$$
Putting the above inequalities together yields
\beaa
P(\max_{1 \leq i \leq n} Y_{5i} \leq 0)
& \leq &
 \left(
 { \la
 \over 2 ( 1-\la_0 +  \sqrt{(1-\la_0)^2 - 4 \rho^2})}
 \right)^{n/2} \\
 & = & \exp\Big\{0.5 n
 \ln \Big( { \la
 \over 2 ( 1-\la_0 +  \sqrt{(1-\la_0)^2 - 4 \rho^2})}
\Big) \Big\} \\
& \leq & \exp( -0.5 n )
\eeaa
(here, $\la= 1.3555\cdots$, $\rho=0.163071\cdots$ and $\la_0= 0.029361\cdots$).
\qed
\section{Lower bound for Theorem \ref{the-1}}\label{sec3}
\setcounter{equation}{0}

Hereafter let $\theta_1(x)=x$, $\theta_2(x)=x^{-1}$, $\theta_3(x)=-x^{-1}$
and $\theta_4(x)=-x$ be the symmetry transformations preserved by
the Gaussian processes $\hat{f}^b_n(x)$ and let
$\oo \gamma_n(x) = \max_{j=1}^4 \gamma_n(\theta_j(x))$
(with the exception of $x=0$ for which $\oo\gamma_n(0)=\gamma_n(0)$).
We begin by noting that, with $I_1=[0,1]$, $I_2=[1,\infty)$,
$I_3=(-\infty,-1]$ and $I_4 = [-1,0]$,
\beqn{3.1}
P(\hat{f}^b_n (x) > \gamma_n(x),\, \forall x \in\reals )
& = & P(\hat{f}^b_n(x) > \gamma_n(x), \; \forall x
\in I_1 \cup I_2 \cup I_3 \cup I_4)
\nonumber \\
& \ge &
\prod_{i=1}^4 P(\hat{f}^b_n (x) > \gamma_n(x), \, \forall x \in I_i)
\nonumber \\
 &\geq& \Bigl[
P(\hat{f}^b_n (x) > \oo\gamma_n(x), \, \forall x \in [0,1])\Bigr]^4
\eeqn
where the first inequality follows by Slepian's lemma
due to the positivity of the covariance $c_n(x,y)$ of $\hat{f}^b_n$,
while the second holds because $c_n(x,y) = c_n(-x,-y) =
c_n(\frac{1}{x}, \frac{1}{y})$. Set
$T=\log n$. The assumptions of Theorem \ref{the-1} imply the
existence of the integers $\log\log T \ll  \tau_T \ll T$ such that
$\delta_n := \sup \{ \oo\gamma_n(x) : x \in [1-\xi_n,1] \} \to 0$ for
$\xi_n = \exp(-\tau_T)$. Recall also our assumption that
$\sup \{ \gamma_n(x) : x \in \reals, n \} \leq M < \infty$.
Applying Slepian's lemma once more yields that,
\beqn{3.2}
&& P(\hat{f}^b_n (x) > \oo\gamma_n(x), \, \forall x \in [0,1])
\nonumber \\
&\ge&
P\Bigl(\inf_{0 \leq x \leq 1-\xi_n} \, \hat{f}^b_n (x) > M \Bigr)
P\Bigl(\inf_{1-\xi_n \leq x \leq 1-n^{-1}} \, \hat{f}^b_n (x) > \delta_n \Bigr)
P\Bigl(\inf_{1-n^{-1} \leq x \leq 1} \, \hat{f}^b_n (x) > M \Bigr)
\nonumber \\
&:=& A_n B_n C_n \;.
\eeqn
Starting with $A_n$, note that for $1 > x \ge y \ge 0$ one has
$$
1\le g (x,y) = \frac{1-xy}{\sqrt{(1-x^2)(1-y^2)}} \le
\sqrt{\frac{1-y}{1-x}}
$$
and hence, by \req{1.5}, taking $x=1-e^{-t}$ and $y=1-e^{-s}$ we see that for
$x,y \in [0,1)$,
\beq{3.0}
c_n (x,y) \ge \sqrt{\frac{1-x\vee y}{1-x \wedge y}} = e^{-|t-s|/2}\,.
\eeq
Recall that $\exp(-|t-s|/2)$ is the covariance of the
stationary Ornstein-Uhlenbeck process
(see \req{o-u}). 
In view of \req{3.0}, we have by Slepian's lemma that
$$
A_n = P\Bigl(\inf_{0\le x \le 1-\xi_n} \hat{f}^b_n(x) > M \Bigr)
\ge P\Bigl(\inf_{0 \le t \le \tau_T} X_t > M \Bigr)\,.
$$
Since $X_t$ is a centered stationary Gaussian process of positive
covariance, yet another application of Slepian's lemma yields that
\beq{3.4}
\liminfn \, \frac{\log A_n}{\log n} \geq
\liminfT \frac{\tau_T}{T} \log P(\inf_{0 \leq t \le 1} X_t > M) = 0\,.
\eeq
(since the random variable $\inf_{0 \le t \le 1} X_t$ is unbounded).

We next turn to the dominant term $B_n$.
Setting $z=1-x, w=1-y$,
for all
$x,y\in [1-\xi_n, 1)$, $n$ large,
it follows from \req{1.5}, \req{1.6} and Lemma~\ref{lem-0} that
$$
c_n(x,y)  \ge \frac{1}{g(x,y)}
 \ge  \frac{2\sqrt{zw}}{z+w} \Bigl[1-
 \frac{(z-w)^2}{\max(z,w)}
\Bigr]
$$
Making yet another change of variables
$z=e^{-t}, w=e^{-s}$, we thus get that
for $\alpha=
e^{-\tau_T}$, in the notations of Lemma~\ref{lem-2.2},
$$
c_n(x,y) \ge \frac{2e^{-\frac{|s-t|}{2}}}{1+e^{-|s-t|}}
\Bigl[1-\alpha (1-e^{-|s-t|})^2\Bigr] = R^{(\alpha)} (s-t) \,.
$$
With $R^{(\alpha)}(0)=c_n(x,x) = 1$, it follows by Slepian's lemma that
$$
B_n \ge P \Bigl( \inf_{0\le t \le T} Y_t^{(\alpha)} > \delta_n \Bigr)
$$
Since $\delta_n \to 0$ and $\alpha_T \log T \to 0$ by our choice of $\tau_T$,
it follows by \req{2.6n} of Lemma~\ref{lem-2.3} that
\beq{3.6}
\liminfn \frac{\log B_n}{\log n } \ge - \frac{b}{4} \;.
\eeq
Finally, we
recall that the sequence
$C_n$ is bounded away from zero by Lemma \ref{lem-2.1}.
Combining \req{3.1}, \req{3.2}, \req{3.4} and \req{3.6}, we thus
arrive at the stated lower bound
$$
\liminfn \frac{\log P^b_{n,\gamma_n}}{\log n} \ge - b \;,
$$
of \req{1.3}.
\qed

\section{Upper bound for Theorem \ref{the-1}}\label{sec4}

\setcounter{equation}{0}

Fixing $\frac{1}{2} > \delta>0$, define the four
disjoint intervals $\cI_1 = [ 1-n^{-\delta}, 1-n^{-(1-\delta)} ]$ and
$\cI_j = \theta_j(\cI_1)$, $j=2,3,4$. Let $V=\bigcup_{j=1}^4 \cI_j$ and
$U=\bigcup_{j=1}^4 \{(x,y): x,y \in \cI_j \}$.

The crucial tool in the proof of the upper bound is the following lemma, whose
proof is deferred to the end of this section:
\begin{lemma}
\label{lem-4.1}
For all $n$ sufficiently large there exist $0 \le \al_n \le n^{-\delta/2}$
such that
\beq{4.1}
c_n (x,y) \le \frac{(1-\al_n)}{g(x,y)} \won_{(x,y)\in U} + \al_n
\quad \quad \forall x,y \in V \,.
\eeq
\end{lemma}
Equipped with Lemma~\ref{lem-4.1}, we show how to complete the proof of the
upper bound.  Let $\{N, b_i^{(j)}, j=1,2,3,4, i=0,\ldots \}$
be independent, identically distributed standard Normal
random variables. For $x \in \cI_1$ consider the infinite
random polynomials
$\hat{f}_\infty^{(j)}(x)= \sqrt{1-x^2} \sum_{i=0}^\infty b_i^{(j)} x^i $
which are for $j=1,2,3,4$ well defined i.i.d. centered
Gaussian processes of covariance function $1/g(x,y)$. Recall that
$g(x,y)$ is invariant to application of each of the invertible transformations
$\theta_j(\cdot)$, $j=2,3,4$ on both $x$ and $y$. Each such transformation
is a one to one map of $\cI_j$ to $\cI_1$. Hence,
the right hand side of \req{4.1} represents
the covariance of the centered Gaussian field
$\widetilde{f}_n(\cdot)$ defined on $V$, of the form
$$
\widetilde{f}_n (x) = \sqrt{1-\al_n} \sum_{j=1}^4
\won_{x\in \cI_j} \hat{f}_\infty (\theta_j(x)) + \sqrt{\al_n} N
$$
Observe that the assumptions of Theorem \ref{the-1} imply that
$\eta_n := n^{-\delta/8} \vee \sup \{ -\gamma_n(x) : x \in V \}$
decay to zero as $n \to \infty$.
With $g(x,x)=1$ for all $x \in V$, relying upon \req{4.1} and
the positivity of $1/g(x,y)$ we get by two applications of Slepian's lemma
that for all $n$ sufficiently large
\beqn{4.2}
P^b_{n,\gamma_n} &=& P\Bigl(\hat{f}^b_n(x) < - \gamma_n(x)
\,, \forall x \in \reals \Bigr) \le
P\Bigl(\sup_{x\in V} \hat{f}^b_n(x) \le \eta_n \Bigr)
\le P\Bigl(\sup_{x\in V} \widetilde{f}_n (x) \le \eta_n \Bigr) \nonumber\\
&\le& P \Bigl( N \le - n^{\delta/8}\Bigr)
 + \prod_{j=1}^4 P\Bigl(\sup_{x\in \cI_j}
 \hat{f}_\infty (\theta_j(x)) \le \frac{2 \eta_n}{\sqrt{1-\alpha_n}}
 \Bigr) \nonumber \\
&\le& e^{-n^{{\delta/4}}/{2}} +
P\Bigl(\sup_{x\in \cI_1} \hat{f}_\infty (x) \le 3 \eta_n\Bigr)^4 \;.
\eeqn
Hence, it is enough to show that
\beq{4.3}
\limsup_{\delta \to 0} \limsupn \frac{1}{\log n} \log
P\Bigl(\sup_{x\in \cI_1} \hat{f}_\infty (x) \le 3 \eta_n \Bigr)
\le - \frac{b}{4} \,.
\eeq
The change of variables
$x=1-z=1-e^{-t}$, $y=1-w=1-e^{-s}$ yields, by \req{1.6} and
Lemma~\ref{lem-0}, that for all sufficiently large $n$ and all $x,y \in \cI_1$,
\beq{4.4}
\frac{1}{g(x,y)} \le \frac{2\sqrt{zw}}{z+w} =
\sech \Bigl(\frac{t-s}{2}\Bigr) \,.
\eeq
For $T=\log n$, $T'=(1-2\delta)T$ and $\ep_{T'}:=3 \eta_n \to 0$,
by yet another application of Slepian's lemma and the stationarity
of the process $Y_t$ of Lemma~\ref{lem-2.2}, it follows from
\req{4.4} that
\beq{ofer-5.5}
P\Bigl(\sup_{x\in \cI_1} \hat{f}_\infty (x) \le 3 \eta_n \Bigr)
\leq P\Bigl(\sup_{t\in [\delta T, (1-\delta)T]} Y_t \le 3 \eta_n \Bigr)
\leq P\Bigl(\sup_{0\le t \le T'} Y_t\le \ep_{T'} \Bigr) \,.
\ee
Consequently, by \req{2.6},
\beq{4.6}
\limsupn \frac{1}{\log n} \log
P\Bigl(\sup_{x\in \cI_1} \hat{f}_\infty (x) \le 3 \eta_n \Bigr)
\le (1-2\delta) \limsupT \frac{1}{T} \log P\Bigl(
\sup_{0\le t \le T} Y_t\le \ep_T \Bigr) = - (1-2\delta) \frac{b}{4}.
\eeq
Taking $\delta \to 0$, we see that \req{4.6} implies \req{4.3}
hence the proof of the upper bound in \req{1.3},
modulo the proof of Lemma~\ref{lem-4.1} which we provide next.

\medskip
\noindent
{\bf Proof of Lemma~\ref{lem-4.1}}
Considering separately $(x,y) \in U$ and $(x,y)\not\in U$,
it is enough by the symmetry relations
$g(x,y) = g(-x,-y) = g(\frac{1}{x}, \frac{1}{y})$ to show that
\begin{subequations}
\label{4.7}
\beqn{4.7a}
&& g(x^n, y^n) \le (1-\al_n) + \al_n g(x,y), \hspace{2cm}
x,y \in \cI_1 \\
\label{4.7b}
&& g(x^n, y^n) \le \al_n g(x,y), \hspace{3cm}
x\in \cI_1, y\in \cI_j, j \neq 1
\eeqn
\end{subequations}
Turning first to \req{4.7a},
recall that $g(x,y)$ is a symmetric function, which equals $1$ on the
diagonal $x=y$. We thus may and shall take without loss
of generality $y>x$. Fixing $x \in \cI_1$, the change of variables
$y=y(\eta)=1-(1-x)(1-\eta)^2$ for $\eta \in (0,1)$ then corresponds to
$\eta=1-\sqrt{w/z}$ where $z=1-x$ and $w=1-y$.
It follows from \req{4.4} that
for all $n$ sufficiently large,
$$
g(x,y) -1 \ge \frac{z+w}{2\sqrt{zw}}-1 =
\frac{1}{2} \Bigl(1-\sqrt{\frac{w}{z}}\Bigr)^2 =
\frac{\eta^2}{2}\,.
$$
Moreover, when $n$ is large enough,
\beq{4.9}
\sqrt{1-x^{2n}} \sqrt{1-y^{2n}} \ge \frac{2}{3}
\quad \quad \forall x,y \in \cI_1
\eeq
So that,
\beq{4.8}
\frac{g(x^n, y^n) -1 }{g(x,y) - 1} \le
\frac{2}{\eta^2}\; \frac{1-x^ny^n- \sqrt{1-x^{2n}} \; \sqrt{1-y^{2n}}}
{\sqrt{1-x^{2n}}\;\sqrt{1-y^{2n}}}
 \le \frac{3}{\eta^2} h(\eta) \,,
\eeq
where
$$
h(\eta):=1-x^n y(\eta)^n -
\sqrt{1-x^{2n}} \; \sqrt{1-y(\eta)^{2n}} \;.
$$
Note that $y(0)=x$, hence $h(0) = 0$. It is not hard to check that
$h'(0)=0$ and
\beaa
h''(\xi) &=&
n y''(\xi) \Bigl[\frac{\sqrt{1-x^{2n}}}{\sqrt{1-y^{2n}}} y^{2n-1}
-x^ny^{n-1} \Bigr] \\
&& +
ny'(\xi)^2 \Bigl[(2n-1) y^{2n-2} \frac{\sqrt{1-x^{2n}}}{\sqrt{1-y^{2n}}}
+ \frac{ny^{4n-2}}{1-y^{2n}} \frac{\sqrt{1-x^{2n}}}{\sqrt{1-y^{2n}}}
- (n-1) x^n y^{n-2} \Bigr] \,,
\eeaa
evaluated at $y=y(\xi)$.
Observing that $|y'(\xi)| \le 2$, $y''(\xi) \in [-2,0]$ and
$x^n\le y(\xi)^n \le y(\eta)^n \le e^{-n^\delta}$ for all $\xi \in [0,\eta]$,
it is easy to check that
there exists a universal finite constant $c_1$ such that
$$
\sup_{\xi \in [0,\eta]} h''(\xi) \le c_1 y(\eta)^n \le c_1 e^{-n^\delta} \,,
$$
for all $n$ large enough and any $x,y \in \cI_1$. Hence,
$h(\eta) \le \frac{1}{2} c_1 e^{-n^\delta} \eta^2$.
Substituting in \req{4.8}, we conclude that
$$
\sup_{x,y\in \cI_1} \frac{g(x^n, y^n) -1}{g(x,y)-1} \le 2 c_1 e^{-n^\delta}
$$
proving \req{4.7a}.

Turning to the proof of \req{4.7b} we assume first that
$x\in \cI_1$ and $y \in \cI_2 \cup \cI_3$. Then, $x,|y|^{-1} \in \cI_1$
with $\sqrt{1-x^2}\sqrt{1-y^{-2}} \leq 1$ and \req{4.9} holding
for $x$ and $y^{-1}$. Moreover,
$x^{n} \vee |y|^{-n} \leq e^{-n^\delta}$, so we have in this case that
$$
\frac{g(x^n, y^n)}{g(x,y)} \leq \frac{3}{2} \frac{|x^n - y^{-n}|}{|x-y^{-1}|}
\leq \frac{3}{2} \sum_{k=0}^{n-1} x^k |y|^{-(n-1-k)} \leq 2 n e^{-n^\delta} \,.
$$
In the remaining case of $x\in \cI_1$ and $y\in \cI_4$ we have that $|y|^n \leq e^{-n^\delta}$,
hence $g(x^n, y^n) \le 2$ while
$$
\frac{1}{g(x,y)} =\frac{\sqrt{1-x^2}\sqrt{1-y^2}}{1+|xy|} \le [ 1-(1-n^{-\delta})^2] \leq
2 n^{-\delta},
$$
thus completing the proof of \req{4.7b}.
\qed

\section{Proof of Theorem~\ref{the-3}}\label{sec7}
\setcounter{equation}{0}

Our proof of Theorem~\ref{the-3} combines the Koml\'os-Major-Tusn\'ady
strong approximation theorem with Theorem~\ref{the-1}.
To this end, note that for every $k$ and $|x| \leq 1$, the sequence
$\{(1-x^2) x^{2j}: j=0,\ldots,k-1\} \cup \{x^{2k}\}$
is a probability distribution,
hence for any real valued $s_j$,
\be\label{able}
\left|s_0 + \sum_{j=1}^k (s_j-s_{j-1}) x^{2j} \right|
= \left|s_k x^{2k} + (1-x^2) \sum_{j=0}^{k-1} s_{j} x^{2j} \right|
\leq \max_{0 \leq j \leq k} |s_j| \,.
\ee
Recall that $E(a_i)=0$ and $E(a_i^2)=1$. Hence,
applying the strong approximation theorem of \cite{kmt} twice
we can redefine $\{a_i, 0 \leq i \leq n-1\}$ on a new probability space
with a sequence of
independent standard normal random variables $\{b_i, \ 0 \leq i \leq n-1\}$
such that for any $p\geq 2$, some $c_p<\infty$, all $t>0$ and $n$,
\be\label{sap}
P\left(\max_{0 \leq j \leq (n-1)/2} \left|\sum_{i=0}^j a_{2i} -\sum_{i=0}^j b_{2i}\right|
\geq t\right) +
P\left(\max_{0 \leq j \leq (n-2)/2} \left|\sum_{i=0}^j a_{2i+1} -\sum_{i=0}^j b_{2i+1}\right|
 \geq t\right) \leq c_p n E|a_0|^p t^{-p} .
\ee
Let
$$
g_k(x) := x^{k-1} f_k(x^{-1}) = \sum_{i=0}^{k-1} a_{i} x^{k-1-i} \;,
$$
and for $k \in \{1,\ldots,n\}$, define
$$
f_k^b(x) = \sum_{i=0}^{k-1} b_{i} x^{i} \;,\qquad \qquad
g_k^b(x) = x^{k-1} f_k^b(x^{-1})\;.
$$
Let $\sigma_k(x):=\sqrt{E(f_k(x)^2)}=
\sqrt{|1-x^{2k}|/|1-x^2|}$, when $|x| \neq 1$ with
$\sigma_k(\pm 1)=\sqrt{E(f_k(\pm 1)^2)}=\sqrt{k}$. Define
$\hat{f}_k(x):=f_k(x)/\sigma_k(x)$ and
$\hat{f}^b_k(x):=f^b_k(x)/\sigma_k(x)$.
As $\sigma_k(x)=\sqrt{E(g_k(x)^2)}$, we shall also use
$\hat{g}_k(x)=g_k(x)/\sigma_k(x)$ and $\hat{g}^b_k(x)=g^b_k(x)/\sigma_k(x)$.
Since
\be
\lbl{able2}
\left|\sum_{i=0}^k a_i x^i -\sum_{i=0}^k b_i x^i\right| \leq
\left|\sum_{j=0}^{\lfloor k/2 \rfloor} (a_{2j}-b_{2j}) x^{2j}\right|
+\left|\sum_{j=0}^{\lfloor (k-1)/2 \rfloor} (a_{2j+1}- b_{2j+1}) 
x^{2j+1}\right|
\ee
we get from (\ref{sap}) by two applications of (\ref{able}) 
(using once $s_j=\sum_{i=0}^j (a_{2i}-b_{2i})$ and once
$s_j=\sum_{i=0}^{
j}(a_{2i+1}-b_{2i+1})$)
that for
all $k \leq n$,
\be\label{sapc1}
P\left(\sup_{|x| \leq 1} \left|f_k(x)-f_k^b(x)\right| \geq  2 t \right) \leq c_p n E|a_0|^p t^{-p}
\ee
The same construction of $\{ b_i \}$ leads by a similar argument also to
\be\label{sapc2}
P\left(\sup_{|x| \leq 1} \left|g_k(x)-g_k^b(x)\right| \geq  4t\right) \leq c_p n E|a_0|^p t^{-p}
\ee
Indeed, bounding $g_k-g_k^b$ amounts to changing $(a_i,b_i)$ to
$(a_{k-1-i},b_{k-1-i})$, resulting with  
using once 
$s_j = \sum_{i=0}^{j} (a_{k-1-2i}-b_{k-1-2i})$ and
once
$s_j = \sum_{i=0}^{j} (a_{k-1-2i-1}-b_{k-1-2i-1})$.
One controls all these as before, but for doubling 
the total approximation error.

In order to apply effectively the strong approximation results,
we need that contributions to the value of $f_n(x)$
come from many variables. This obviously is easier for
$||x|-1|$ 
small. In order to avoid appearance
of zeros in other locations, we decompose $f_n$ to the dominant
``bulk term", which
will not be too negative {\it everywhere} and will be rather far from 
$0$ for 
$||x|-1|$ 
small, and 
to ``boundary terms", which involve a small number of coefficients
and thus can be made to have prescribed 
positive values with a not too small probability.

In order to define precisely the different regions considered for
values of $x$ and the splitting into bulk and boundary terms, 
 we introduce, for $n$ large enough odd integers,
 a few $n$-dependent parameters as follows:
\beq{rhodef}
\begin{array}{lll}
p_n:&p_n\uparrow \infty, c_{p_n} E|a_0|^{p_n}\leq n&
c_p \,\mbox{\rm is the KMT constant in \req{sap}.}\\
\vp=\vp_n:&\vp_n\downarrow 0, \vp_n\geq 
\max\{20/p_n,(\log n)^{-1/2}\},&\vp_n\,\, \mbox{\rm is taken as the
 smallest possible}\\
& 2n^{3\vp_n}=2^j\,\, \mbox{\rm for some $j$ integer}&
 \mbox{\rm 
  value satisfying constraints.
}\\
m=m_n:&m_n\to\infty, m_n=2n^{3\vp_n}& m_n \,\,
\mbox{\rm is an integer power of $2$.}\\
\bar{\gamma}_n(x):&
\bar{\gamma}_n\to 0, 
\bar{\gamma}_n(x)=\max \{0,\gamma_n(x),\gamma_n(x^{-1})\}&
\gamma_n(x) \,\,
\mbox{\rm as in statement of theorem.}\\
\rho_n:&\rho_n\to 0,
\rho_n = \sup_{|x| \leq 1-m^{-1}} \; \{ \sigma_n(x) \bar{\gamma}_n(x) \}&
 \rho_n\leq c n^{-\delta/2},\,\, \mbox{\rm some 
finite $c>0$}.\\ 
r=r_n:&c n^{-\delta/2}& \mbox{\rm for $n\geq 3m$; $c$ is as in bound on 
$\rho_n$.}\\
\xi_n(x):&\xi_n(x)=6x^m \sigma_{n-2m}(x)\bar{\gamma}_n(x)&
\end{array}
\ee
In order to state the decomposition alluded to above,
first 
 partition the interval $[-1,1]$ to
$\cI = \{ x : |x| \geq 1-0.5 n^{- \vp} \}$ and $\cI^c = [-1,1]
\setminus \cI$.
We note that
$2r+\xi_n(x) \geq \sigma_n(x) \bar{\gamma}_n(x)$ for all $x \in \cI$.
Next, let
$f_n=f_n^L+f_n^M+ f_n^H$ where
\beq{LMH}
f_n^L(x)=\sum_{i=0}^{m-1} a_i x^i, \qquad
f_n^M(x)=\sum_{i=m}^{n-1-m} a_i x^i, \qquad
f_n^H(x)=\sum_{i=n-m}^{n-1} a_i x^i \;.
\eeq
Similarly, we let $g_n=g_n^L+g_n^M+ g_n^H$ with
$g_n^L (x)=x^{n-1} f_n^L (x^{-1})$, etc.
With these definitions, we have the inclusions
\beaa
\{ \hat{f}_n(x) > \gamma_n(x), \quad \forall \ x \in \reals\}
&\supset& \{ \hat{f}_n(x) >\bar{\gamma}_n(x), \; \hat{g}_n(x) > \bar{\gamma}_n(x),
\quad \forall \ x \in [-1,1]\} \\
&\supset & \{ f_n^M(x) > \xi_n(x),\; g_n^M(x) > \xi_n(x),
\quad \forall \ x \in \cI \} \\
&& \cap \{ f_n^M(x) >-r,\; g_n^M(x) >-r, \quad \forall \ x \in \cI^c \}\\
&& \cap \{ f_n^L(x) >3r, \; g_n^L(x) \geq -r, \quad \forall  x \in [-1,1]\}\\
&& \cap \{ f_n^H(x) \geq -r, \; g_n^H(x) > 3r, \quad \forall x \in [-1,1]\}
\eeaa
($f_n^M, g_n^M$ are the ``bulk terms'' whereas $f_n^L,g_n^L,f_n^H,g_n^H$
are the ``boundary terms'').
Since the polynomial pairs $(f_n^L,g_n^L)$, $(f_n^M,g_n^M)$ and
$(f_n^H,g_n^H)$ are mutually independent, it follows that
\beqn{breakup}
P_{n,\gamma_n} &=&
P\Big( \hat{f}_n(x) >\gamma_n(x), \quad \forall \ x \in \reals \Big) \nn \\
& \geq & P\Big(
 \{f_n^M (x) >\xi_n(x), \; g_n^M(x) >\xi_n(x),
\;\, \forall \ x \in \cI\}\cap
\{f_n^M(x) >-r, \; g_n^M(x) > -r, \;\, \forall \  x \in \cI^c\} \Big) \nn \\
&& P\Big( f_n^L(x) >3r, \;\; g_n^L(x) \geq -r, \quad
\forall \ x \in [-1,1] \Big) \nn \\
&& P\Big( f_n^H(x) \geq -r, \;\; g_n^H(x) > 3r, \quad
\forall \ x \in [-1,1] \Big)
\eeqn
Note that $g_n^M$ and $f_n^M$ are identically distributed, as are the
polynomial pairs $x^{-m} (f_n^M,g_n^M)$ and $(f_{n-2m}, g_{n-2m})$. Thus,
we have that
\beaa
\lefteqn{
P\Big( 
\{f_n^M (x) >\xi_n(x), \;\; g_n^M(x) >\xi_n(x),
\quad \forall \ x \in \cI\}\cap\{
f_n^M(x) >-r, \;\; g_n^M(x) > -r, \quad \forall \ x \in \cI^c\} \Big) } \\
&\geq&
P\Big( \hat{f}_{n-2m} (x) > 6 \bar{\gamma}_n(x),
\; \hat{g}_{n-2m} (x) >6 \bar{\gamma}_n(x), \;\; \forall \ x \in \cI \Big) -
 2 P\Big(f_n^M(x) \leq -r , \; \mbox{for some} \ x \in \cI^c \Big)
\\
&:= & Q_1 - 2 Q_2
\eeaa
Since the polynomial pairs $(f_n^L,g_n^L)$, $(g_n^H,f_n^H)$ and
$(f_m,x^{n-m} g_m)$ have identical laws, it now follows that
$$
P_{n,\gamma_n} \geq (Q_1-2 Q_2) (Q_3-Q_4)^2\,,
$$
where
$$
Q_3 := P\Big( f_m(x) >3r, \;\; \forall \ x \in [-1,1], \quad
x^{n-m} g_m(x) \geq -r, \;\; \forall \ |x| \in [1 - m^{-1},1] \Big)\,,
$$
and 
$$ 
Q_4 :=
P\Big( x^{n-m} g_m(x) \leq -r, \; \mbox{for some} \ |x| \leq 1 - m^{-1} \Big).
$$
To deal with the dominant term $Q_1$, we
consider \req{sapc1} and \req{sapc2} for $p=p_n$ as above,
$k=n-2m$, and $t=n^{\vp/4}$. Noting that
$\eta_n=\sup \{ 6\bar{\gamma}_n(x) + 4 t/\sigma_{n-2m}(x) : x \in \cI\}$
approaches zero as $n \to \infty$, we get that for all $n$ large enough,
\bea
\label{q1bd}
Q_1 &\geq&
P\Big( \hat{f}^b_{n-2m} (x) > \eta_n, \quad  \hat{g}_{n-2m}^b(x) >  \eta_n,
\quad \forall  \ x \in \cI \Big) - 2 n^{-3} \nn \\
&\geq& P\Big( f_{n-2m}^b (x) > \eta_n \sigma_{n-2m}(x), \quad
\quad \forall \ | |x| - 1 | \leq n^{-\vp_n} \Big) - 2 n^{-3}
\nn \\ & \geq & (n-2m)^{-b+o(1)},
\eea
where the last inequality follows
by applying Theorem \ref{the-1} for threshold
$\eta_n \to 0$ for $| |x| - 1 | \leq n^{-\vp_n}$ and zero otherwise.

Turning to estimate $Q_2$, recall that $f_n^M$ has the same
distribution as $x^m f_{n-2m}$ and
$m=2 n^{3\vp} $.
Recall also that $\vp \geq (\log n)^{-1/2}$, implying that
$n^c \exp(-n^\vp) \to 0$ for any fixed $c<\infty$. Hence,
for all $n$ large enough,
\beaa
Q_2 &\leq& P\Big( \sup_{x  \in \cI^c} |x|^m |f_{n-2m}(x)| \geq r \Big)
\\
&\leq &
P\Big( \sup_{x \in \cI^c} |f_{n-2m}(x)|\geq  cn^{-\delta/2}\Big(1-
\frac{1}{2 n^\vp}
\Big)^{2n^{3\vp}}\Big)
\leq P\Big( \sup_{x \in \cI^c} |f_{n-2m}(x)|\geq  2 \exp(n^{\vp})
\Big).
\eeaa
Observe that for any $x, y \in [-1,1]$,
$$
E\left(\left|f_{n-2m}(x)-f_{n-2m}(y)\right|^2\right)
\leq \sum_{i=1}^{n} (x^i -y^i)^2 \leq (x-y)^2 n^3 \;.
$$
Recall the following well known lemma (see \cite{stout} for a proof).
\begin{lemma}\label{moment}
Let $\{T_x, x \in [a, b]\}$ be an a.s. continuous
stochastic process with $T_a=0$.
Assume that
$$ \forall \ \ x, y \in [a,b], \quad E|T_x-T_y|^2 \leq K (x-y)^2.
$$
Then, we have
$$
E\Big( \sup_{x \in [a,b]} T_x^2 \Big) \leq 4 K (b-a)^2.
$$
\end{lemma}
Applying Lemma \ref{moment} for $T_x=f_{n-2m}(x)-f_{n-2m}(0)$,
first when $x \in [0,1]$, then when $x \in [-1,0]$, we get by
Markov's inequality that for all $n$ large enough,
and any 
$c_1<\infty$
(for our use, $c_1=3$ will do),
\be
Q_2  \leq P \Big(|a_0| \geq \exp(n^{\vp}) \Big) +
 P\Big( \sup_{|x| \leq 1} T_x^2 \geq \exp(2 n^{\vp}) \Big)
\leq \exp(-2 n^{\vp}) (1+E( \sup_{|x| \leq 1} T_x^2 )) = o(n^{-c_1}) \,.
\label{i2}
\ee
Recall that $m=2 n^{\vp_n} $ and $\vp_n \to 0$, so
with $g_m$ and $f_m$ of identical law, it follows that
for all $n$ large enough,
$$
Q_4 \leq P\Big( \sup_{|x| \leq 1 - m^{-1}} |x|^{n-m} |f_m(x)| \geq r \Big)
\leq P\Big( \sup_{|x| \leq 1 - m^{-1}} |f_m(x)|\geq  2 \exp(\sqrt{n}) \Big).
$$
Similarly to the derivation of (\ref{i2}), by
twice applying Lemma \ref{moment} for $T_x=f_{m}(x)-f_{m}(0)$,
then using Markov's inequality, we get that $Q_4 \leq \exp(-n^{1/3})$
for all $n$ large enough. The lower bound $P_{n,\gamma_n} \geq n^{-b+o(1)}$
in Theorem \ref{the-3}
is thus a direct consequence of the bounds \req{q1bd}, \req{i2}
and Lemma \ref{lower} below which provides the estimate $Q_3\geq n^{-c_2 
\vp}$ with 
$\vp=\vp_n \to 0$ and $c_2<\infty$ fixed.

Turning to the upper bound $P_{n,\gamma_n} \leq n^{-b+o(1)}$ in
Theorem \ref{the-3}, let $\eta_n:=\inf\{ \gamma_n(x) : ||x|-1|\leq n^{-\vp}
 \}$. Recall that $\eta_n \to 0$ by our assumptions. Then,
similarly to the derivation of \req{q1bd}, now with $m=0$, we see that
for all $n$ large enough
\beaa
P_{n,\gamma_n} &=&
P\Big( \hat{f}_n(x) >\gamma_n(x), \quad \forall \ x \in \reals \Big) \\
&\leq& P\Big(
\hat{f}_n(x) >\eta_n,
\quad \hat{g}_n(x) >\eta_n, \quad \forall \ x \in \cI )  \\
&\leq&
P\Big( \hat{f}^b_n(x) > \eta_n - n^{-\vp/8},
\quad \hat{g}_n^b(x) >\eta_n  -n^{-\vp/8},
\quad \forall \ x \in \cI ) + 2 n^{-3} \\ &\leq& n^{-b+o(1)}
\eeaa
(the last inequality follows by Theorem \ref{the-1} for a
threshold $\eta_n-n^{-\vp_n/8}$ when $x \in \cI \cup \{x: x^{-1} \in \cI \}$
and $-\infty$ otherwise).
\qed

\begin{lemma} \label{lower}
Suppose $a_i$ are i.i.d. with  $E(a_0)=0$ and $E(a_0^2)=1$.
There exists $c<\infty$ such that for all $m=2^{k+1}$ and $k$ large enough,
\be
P\Big( f_m(x) > m^{-2}, \;\; \forall \ x \in [-1,1]\;, \quad \quad
x g_m(x) \geq 0, \;\;  \forall \ |x| \in [1-2^{-k},1] \Big)\,
\geq m^{-c}
\label{i3}
\ee
\end{lemma}

\proof 
Define the intervals
$\cJ_{j}=\{ x : 1-2^{-j} \leq |x| \leq 1-2^{-j-1} \}$ for $j=1,\ldots,k-1$
and $\cJ_{k}=\{ x : 1-2^{-k} \leq |x| \leq 1 \}$. Throughout this proof,
$l_j:=2^j$ for integer $j$, and
complements are taken inside the interval $[-1,1]$.

The proof of the lemma is based on decomposing $f_m$ to a sum (over
a number of terms logarithmic in $m$) of
polynomials 
$f^j$, such that for each $x\in \cJ_j$,
$f^j(x)$ is large while $f^i(x), i\neq j$ are not too large; at the
same time,  
$g_m(x)$ is decomposed to a sum of polynomials  
all but the highest order of which are large and positive on $\cJ_k$, 
while the latter is not too negative on $\cJ_k$. 
Unfortunately, we need to introduce a few constants
in order to define explicitly this decomposition.

Note first that
for some $c_0<\infty$ which does not depend on $k$, 
\be
\label{v22}
(c_0-1) 2^{j/2} x^{l_j} -
\sum_{i=4,i \neq j}^{k} 2^{i/2} x^{l_i} \geq 0\; , \quad
\quad \forall x \in \cJ_j, \;\; j=4,\ldots,k
\ee
Define $c_1=c_0+1$. 
In Lemma \ref{ofer-mar} below, we define a constant
$\theta_1=\theta_1(c_1)>0$. Define then $\theta=P(|N|\leq 1)\theta_1/2>0$
where $N$ is a standard Gaussian random variable of 
zero mean and unit variance. 
Since $E(a_0)=0$, $E(a_0^2)=1$, we can use
Strassen's weak approximation theorem
(see \cite{Str64} or \cite[Page 89]{cs}), to
deduce the existence of independent standard normal
random variables $\{b_i, i \geq 0 \}$ such that, for all $j\geq j_0$,
\be
P\left(
\max_{0 \leq \ell \leq  2^j} \left|\sum_{i=0}^{\ell} a_{2i}-
\sum_{i=0}^{\ell} b_{2i}\right|
+\max_{0 \leq {\ell} \leq  2^j}\left|\sum_{i=0}^{\ell} a_{2i+1}-
\sum_{i=0}^\ell b_{2i+1} \right| \geq 2^{j/2-3} \right)\leq \theta.
\lbl{qn3}
\ee
Finally, since $Ea=0$ and $Ea^2=1$ there exists $\al>0$ such
that $P(|a-\alpha| \leq \delta)>0$. Fixing such $\alpha$,
define $s>j_0$
such that 
\be
\label{v3}
{\al \over 10} -
\sum_{i=s}^\infty 2^{i/2} x^{l_i} \geq 0, \quad \quad \forall x \in \cJ_0
:=\{ x : |x| \leq 1-2^{-s} \}.
\ee
Such an $s$ always exists because the sum in \req{v3} tends to $0$
in $s$. Note that $s$ does not depend on $k$ and all estimates
above are valid uniformly for all $k$ large enough.
We write $l:=l_s$ and note that $\{\cJ_0,\cJ_s,\cJ_{s+1},
\ldots,\cJ_k\}$ form a partition of the interval $[-1,1]$.
We keep $s$ fixed throughout the rest of the proof.

As mentioned above, 
the proof of the lemma is based on decomposing $f_m$ to a sum (over
$k-s+2$ 
terms, i.e. a number of terms logarithmic in $m$) of
polynomials $f^j, j=0,s,s+1,\ldots,k$, 
while 
decomposing $g_m(x)$ to a 
similar sum of $k-s+2$ polynomials.
Specifically, we
write
$$
f_m(x) =f^0 (x) + \sum_{j=s}^{k} x^{l_j} f^j(x)
$$
where $f^0 = f_{l_s}$ and
$$
f^j(x) = \sum_{i=0}^{2^j-1} a_{i+l_j} x^i \;.
$$
Similarly, 
$$
g_m(x) =x^{m-l_s} g^0 (x) + \sum_{j=s}^{k} x^{m-l_{j+1}} g^j(x)
$$
where $g^0 = g_{l_s}$ and
$$
g^j(x) = \sum_{i=0}^{2^j-1} a_{i+l_j} x^{2^j-1-i} \;.
$$

One checks that  for $k$ large enough,
it holds that
\be
\label{v1}
m^{-2} \leq
\min\{{\al \over 10},\inf_{x \in \cJ_j, j=s,\ldots,k}\, 2^{j/2} x^{l_j}\}\,.
\ee
Moreover, 
by \req{v22},
 for 
all
$k \geq s $, 
\be
\label{v2}
(c_0-1) 2^{j/2} x^{l_j} -
\sum_{i=s,i \neq j}^{k} 2^{i/2} x^{l_i} \geq 0\; , \quad
\quad \forall x \in \cJ_j, \;\; j=s,\ldots,k
\ee
It follows that
\bea
\label{dec1}
\{ f_m(x) > m^{-2}, \;\; \forall \ x \in [-1,1] \} & \supset &
\bigcap_{j=s}^k
\{ f^j(x) > c_0 2^{j/2}, \; \forall \ x \in \cJ_j, \;\,
 f^j(x) \geq -2^{j/2}, \; \forall \ x \in \cJ^c_j \} \nn \\
&\bigcap &
\{ f^0(x) > {\al \over 5}, \;\; \forall \ x \in \cJ_0, \quad
 f^0(x) \geq 0, \;\; \forall \ x \in \cJ^c_0 \}\,.
\eea
Note that for all $x \in [-1,1]$,
\be
\label{dec2}
\{ x g_m(x) \geq 0 \} \supset  \{ x g^0(x) \geq -2^{k/2} \}
\bigcap_{j=s}^{k} \{ x g^j(x) \geq c_0 2^{j/2} \} \;.
\ee
The polynomial pairs $(f^j,g^j)$, $j=0,s,\ldots,k$ are mutually
independent, with $(f^0,g^0)$ having the same law as
$(f_l,g_l)$, while $(f^j,g^j)$ has
the same law as $(f_{2^j},g_{2^j})$ for each $j \neq 0$.
It thus follows from  (\ref{dec1}) and (\ref{dec2}) that
\bea
\label{qn1}
&& P\Big( f_m(x) > m^{-2}, \; \forall \ x \in [-1,1] \,,\quad\quad
x g_m(x) \geq 0, \;  \forall \ x \in \cJ_k \Big) \nn \\
&\geq&
P\Big( f_{l}(x) > {\al \over 5}, \;\; \forall \ x \in \cJ_0, \quad
 f_l(x) \geq 0, \;\; \forall \ x \in \cJ^c_0, \quad
 x g_l(x) \geq -2^{k/2}, \;\; \forall \ x \in \cJ^k \Big) \nn \\
&\prod_{j=s}^{k}&
P\Big( f_{2^j} (x) > c_0 2^{j/2}, \;\, \forall \ x \in \cJ_j, \quad
 f_{2^j}(x) \geq -2^{j/2}, \;\; \forall \ x \in \cJ^c_j, \quad
 x g_{2^j}(x) \geq  c_0 2^{j/2}, \;\, \forall \ x \in \cJ_k \Big) \nn \\
&:=& 
\eta_{s,k} \prod_{j=s}^{k} q_j \,,
\eea
where
\be
\label{q00mar}
\eta_{s,k} := P\Big( f_{l}(x) > \al/5, \;\; \forall \ x \in \cJ_0, \quad
 f_l(x) \geq 0, \;\; \forall \ x \in \cJ^c_0, \quad
 x g_l(x) \geq -2^{k/2}, \;\; \forall \ x \in \cJ^k \Big)\,. 
\ee
We first show that $\eta_{s,k}$ is uniformly (in $k$, for $k\geq 2\log_2(2\al l)$),
bounded away from zero, and then provide a uniform (in $k$)
bound (independent of $j$)
on $q_j$. Toward the first goal,
let $Q_s(x):=\al (1+x+\cdots+x^{l_s-1})$, noting that 
$Q_s(x)$ is monotone increasing on $[-1,1]$, with
$Q_s(-(1-2^{-s}))\geq \al/4$, implying
that $Q_s(x) \geq \al/4$ for all $x \in \cJ_0$. Thus, for each $s \geq 1$
there exists $\delta_s \in (0,\al)$ such that $f_{l_s}(x) > \al/5$ whenever 
$x \in \cJ_0$ and $|a_i-\al| \leq \delta_s$ for $i=0,\ldots,l_s-1$.
Further taking
$a_{2i} \geq a_{2i+1} \geq 0$ for $i=0,\ldots,2^{s-1}-1$ guarantees
that $f_{l_s}(x) \geq 0$
for all $x \in [-1,1]$. 
Considering only such $\{a_i\}$, we also have that
$ |xg_l(x)|\leq 2\al l$, and hence, combining the above and using
$2^{k/2}\geq 2\al l$, we have that 
\be
\label{q0}
\liminf_{k \to \infty} \eta_{s,k} >0\,.
\ee 

To estimate $q_j$, we note that
$x g^b_{2^j}=x^2 g^b_{2^j-1}+x b_{2^j-1}$. Thus, combining 
\req{able}, \req{able2} and \req{qn3},
it follows that for all
$j \in \{s,\ldots,k\}$, 
\be
\lbl{qjbd}
q_j \geq P(|b_{2^j-1}| \leq 1) q_j^{b} - \theta
\ee 
where 
$$
q_j^b := P\Big( f^b_{2^j-1} (x) \geq c_1 2^{j/2}, \; \forall
\ x \in \cJ_j, \;\;
f^b_{2^j-1}(x) \geq -2^{(j-1)/2}, \; \forall \ x \in \cJ^c_j, \;\;  
g^b_{2^j-1}(x) \geq  c_1 2^{j/2}, \; \forall \ x \in \cJ_k \Big)\,.
$$
for say, $c_1=c_0+1$.
Slepian's lemma thus yields that for all $k \geq j \geq s$,
$$
q_j^b \geq P\Big( f^b_{2^j-1} (x) \geq c_1 2^{j/2},
\; \forall x \in [1-2^{-j},1] \Big)^4
  P\Big( f^b_{2^j-1} (x) \geq -2^{(j-1)/2},
\; \forall  x \in [0,1-2^{-j}] \Big)^2:=\bar q_j^b.
$$
Note that $\bar q_j^b$ does not depend on $k$, and in fact it depends
on $c_1$ and $j$ only. The following lemma provides estimates on 
$\bar q_j^b$ while defining the constant $\theta_1$:
\begin{lemma}
\label{ofer-mar}
There exists a constant $\theta_1>0$ such that for
all $j\geq 4$,
$$\bar q_j^b\geq \theta_1\,.$$
\end{lemma}
Applying \req{qjbd} using
$\theta = {1 \over 2} P(|b|\leq 1) \theta_1$
then leads
to $q^j \geq \theta$ for
all $j \geq s$.  In view of
\req{q0} and \req{qn1} this proves \req{i3}.
\qed

\noindent
{\bf Proof of Lemma \ref{ofer-mar}:}
Note that $\sigma_{2^j-1}(x) \geq 2^{j/2-1}$
when $x \in [1-2^{-j},1]$, hence
by Lemma \ref{lem-2.1}, for some $\xi_1>0$ and all $j$ large enough,
$$
P\Big( f^b_{2^j-1} (x) \geq c_1 2^{j/2}, \;\, \forall x
\in [1-2^{-j},1] \Big)
\geq 
P\Big( \hat{f}^b_{2^j-1} (x) \geq 2 c_1, \;\, \forall x \in
[1-1/(2^j-1),1] \Big)
\geq \xi_1 \;.
$$
Note that $\sigma_{2^j-1}(x) \leq 1/\sqrt{1-x} \leq 2^{i/2}$
for $x \in [0,1-2^{-i}]$.
Hence, by Slepian's lemma and \req{3.0}, we have that for the 
Ornstein-Uhlenbeck process $X_t$ of \req{o-u},
\beaa 
\lefteqn{
P\Big( f^b_{2^j-1} (x) \geq -2^{(j-1)/2}, \forall  x \in [0,1-2^{-j}] \Big)}
\\
& \geq & \prod_{i=0}^{j-1} P\Big( f^b_{2^j-1} (x) \geq -2^{(j-1)/2}, \,
\forall x \in [1-2^{-i},1-2^{-i-1}] \Big) \\
& \geq &
\prod_{i=0}^{j-1}
P\Big( X_t  \geq -2^{(j-i)/2-1},  \; \forall t \in [i\ln 2, (i+1)\ln 2]\Big) \\
& = &
\prod_{l=1}^j
P\Big( \inf_{0 \leq t \leq \ln 2}\, X_t \geq -2^{l/2-1}\Big) \\
&\geq &\prod_{l=1}^\infty
\Big[ 1-  P(\sup_{0 \leq t \leq \ln 2}\, X_t \geq 2^{l/2-1}) \Big] = : \xi_2 
\eeaa
and $\xi_2>0$ since $E(\sup_{t\in [0,\ln 2]} X_t) < \infty$.
This completes the proof.
\qed

\section{Proof of Theorem \ref{theo-1int}}
\label{sec4new}
\setcounter{equation}{0}

Part a) of Theorem \ref{theo-1int}
is a direct consequence of Theorem \ref{the-3} with $\gamma_n=0$.
Thus, it only remains to prove part b).
Fixing $\mu \neq 0$, it is easy to see that
$$
P_n^\mu=P(\hat{f}_n(x) \neq -\mu \kappa_n(x), \quad \forall x \in \reals)\,,
$$
where the nonrandom
$\kappa_n(x)=(\sum_{i=0}^{n-1} x^i)(\sum_{i=0}^{n-1} x^{2i})^{-1/2}$ are
strictly
positive and
$\hat{f}_n(x)$ are the normalized polynomials that correspond to 
$a_i$ of zero mean. With 
$\widetilde{P}_n$ for the value of $P_n$ when coefficients $\{-a_i\}$
are used instead of $\{a_i\}$, it is easy to see that
$$
P_n^\mu = P_{n,-\mu \kappa_n} + \widetilde{P}_{n,\mu \kappa_n} \,.
$$
Consequently, we may and shall assume without loss of generality that $\mu>0$,
proving only that $P_{n,-\mu \kappa_n} = n^{-b/2 + o(1)}$.
Observe that $\kappa_n(1)=\sqrt{n}$, $\kappa_n(-1)=1/\sqrt{n}$ and
if $|x| \neq 1$ then,
\beq{7.1}
\kappa_n(x)=\kappa_n(x^{-1})=\kappa_n(-x)^{-1}
=\left[\frac{(1+x)(1-x^n)}{(1-x)(1+x^n)}\right]^{1/2} \;.
\eeq
Moreover, there exists $c=c(\mu)>0$ such that for all $n$ large enough,
\beq{7.3}
\mu \kappa_n(x) \geq n^{-\vp/8}+ \frac{c}{\sqrt{1-x}} \quad \quad \forall
x \in [0,1-n^{-1}] \,.
\eeq
For an upper bound on $P_{n,-\mu \kappa_n}$
let $\cI_- = [-1,-1+0.5 n^{-\vp}]$
be the subset of $\cI$ of Section \ref{sec7} near the point $-1$ and
$V_-=\cI_3 \cup \cI_4$ be the (corresponding) subset of $V$ of
Section \ref{sec4}. It is easy to check that
$$
\sup\{ \kappa_n(x) : x \in \cI_- \} \leq c_1 n^{-\vp/2}
$$
for some $c_1<\infty$ and all $n$. Hence, applying
the arguments of Section \ref{sec7} followed by those of Section \ref{sec4}
with $\cI_-$ replacing $\cI$ and $V_-$ replacing $V$, respectively,
results in the stated upper
bound $P_{n,-\mu \kappa_n} \leq n^{-b/2+o(1)}$.
Turning to prove the corresponding
lower bound on $P_{n,-\mu \kappa_n}$, let $\cI_+=[1-0.5 n^{-\vp},1]$
denote the subset of $\cI$ near the point $+1$. 
It follows from \req{7.1} that $\rho_n$ of (\ref{rhodef}) is zero for
$\gamma_n=-\mu \kappa_n$, allowing for the use of $r_n=n^{-1}$ and
$\xi_n(x)=-\mu x^m \sigma_{n-2m}(x) \kappa_n (x) \leq 0$ in (\ref{breakup}).
We then deal with
the terms $Q_2$, $Q_3$ and $Q_4$ as in Section \ref{sec7}. For the dominant
term $Q_1$, instead of (\ref{q1bd}) we have in view of (\ref{7.3}) that 
\beqn{q1nbd}
Q_1 &\geq&
P\Big( \hat{f}^b_{n'} (x) > n^{-\vp/8}, \quad  \hat{g}_{n'}^b(x) >  n^{-\vp/8},
\quad \forall  \ x \in [-1,0] \cup [1-\frac{1}{n'},1], \nn \\
&& \quad \;\; \hat{f}^b_{n'} (x) > -\frac{c}{\sqrt{1-x}}, \quad  
 \hat{g}^b_{n'} (x) > -\frac{c}{\sqrt{1-x}}, 
\quad \forall  \ x \in [0,1-\frac{1}{n'}] \Big) - 2 n^{-3} \nn \\ 
& = & \widetilde{Q}_1 - 2 n^{-3},
\eeqn
where $n':=n-2m = n (1+o(1))$. By Slepian's lemma, similarly to \req{3.1} we
see that
\beqn{7.2}
\widetilde{Q}_1 &\geq&
P\Big( \hat{f}^b_{n'} (x) >
-\frac{c}{\sqrt{1-x}} \quad \forall  \ x \in [0,1-\frac{1}{n'}] \Big)^2
\nn \\ 
&&
P\Big( \hat{f}^b_{n'} (x) > n^{-\vp/8} \quad \forall  \ x \in [0,1] \Big)^2 
P\Big( \hat{f}^b_{n'} (x) > 1 \quad \forall  \ x \in [1-\frac{1}{n'},1] \Big)^2
\nn \\
&:=& (A_{n'})^2 (B_{n'})^2 (C_{n'})^2 \,.
\eeqn
The sequence $C_{n'}$ is bounded away from zero by Lemma \ref{lem-2.1}.
Moreover,
it is shown in Section \ref{sec3} that $B_n \geq n^{-b/4+o(1)}$. 
In view of \req{q1nbd} and \req{7.2}, it thus suffices to show that 
the sequence $A_n$ is bounded below by some $A_\infty>0$ in order to
conclude that
$P_n^\mu \geq P_{n,-\mu \kappa_n} \geq n^{-b/2+o(1)}$ and complete the proof of
part b) of Theorem \ref{theo-1int}.
To this end, recall that the function
$\sqrt{(1-x\vee y)/(1-x \wedge y)}$
in the right side of \req{3.0} is the covariance of the process
$W_{1-x}/\sqrt{1-x}$. Consequently, we have by \req{3.0} and
Slepian's lemma that
$$
A_n \geq P(W_{1-x} > - c, \, \forall x \in [0,1-n^{-1}])
\geq P( \inf_{0 \leq x \leq 1} W_x > - c) = A_\infty >0 \,,
$$
as needed.
\qed

\section{Upper bound for Theorem \ref{theo-2}}\label{sec6}

\setcounter{equation}{0}

Fixing small $\delta>0$ and integers $k_n = o(\log n/\log \log n)$,
it suffices for the upper bound in Theorem \ref{theo-2}
to provide an $n^{-(1-2\delta)b+o(1)}$ upper bound on the probability $q_{n,k}$
that $f(x)=\hat{f}_{n+1}(x)$
has at most $k=k_n$ zeros in the set $V=\cup_{i=1}^4 \cI_i$ of
Section \ref{sec4}. To this end, let
$x=\theta_i(1-e^{-t})$ within $\cI_i$, $i=1,2,3,4$,
where $\theta_i$ is defined in Section \ref{sec3}.
With $T =\log n$ cut the range $t \in [\delta T, (1-\delta) T]$
for each $\cI_i$
to $(1-2\delta)T$ unit length intervals, denoting by
$J_{(i-1)T+1}, \ldots, J_{iT}$ the corresponding image in $\cI_i$.
If $f(x)$ has $\ell$ zeros in some $\cI_i$,
then there must exist $j_1,\ldots,j_\ell$
such that $f(x)$ has a constant sign $\ss \in \{-1,1\}$ on each of the ``long''
subintervals obtained by deleting $J_{(i-1)T+j_1}, \ldots, J_{(i-1)T+j_\ell}$
from $\cI_i$. We partition the event that
$f(x)$ has at most $k$ zeros in $V$ according to the possible
vector ${\bj}=(j_1,\ldots,j_k)$ of ``crossing indices''
among the $4T$ intervals $\{J_1,\ldots,J_{4T}\}$ and the
possible signs $\ss_m \in \{-1,1\}$ of $f(x)$ on the resulting
long subintervals $L_m$, $m=1,\ldots,k+4$ within $V$. Let
$$
q_{n,\bs,\bj} =
P\Bigl( \min_{m=1}^{k+4} \inf_{x \in L_m} \ss_m f(x) > 0\Bigr)\,,
$$
for $\bs= (\ss_1,\ldots,\ss_{k+4})$. Since
$$
q_{n,k} \leq \sum_{\bj} \sum_{\bs} q_{n,\bs,\bj} \,,
$$
and the number of choices of $\bj$ and $\bs$ is at
most $2^{k+4}(4T)^k=n^{o(1)}$, it suffices to show that
\beq{6.1}
\max_{\bs,\bj}  q_{n,\bs,\bj} \leq
n^{-(1-2\delta)b+o(1)}\,.
\eeq
Applying the coupling of $(f_k,g_k)$ and $(f_k^b,g_k^b)$ as provided
in (\ref{sapc1})
and (\ref{sapc2}) for $t=n^{\delta/4}$ and $p=16/\delta$, 
we see that for all ${\bf j}$ and $\bs$,
\beq{6.1new}
q_{n,\bs,{\bf j}} \leq 
P\Bigl( \min_{m=1}^{k+4} \inf_{x \in L_m} \ss_m \hat{f}_{n+1}^b (x)
> -n^{-\delta/8}\Bigr)
+c n^{-3} := q^b_{n,\bs,{\bf j}} +c n^{-3} \,,
\eeq
where $c<\infty$ depends on $\delta$ but
is independent of ${\bf j}$, $\bs$ and $n$. 
Thus, the proof reduces to the Gaussian case.

Suppose first that $n$ is even. 
The covariance function of $\ss_m \hat{f}^b_{n+1}(x)$ on $V' = \cup_m L_m$
is then
given by $\ss_m \ss_l c_{n+1}(x,y)$ for $x \in L_m$ and $y \in L_l$. Since
$c_{n+1}(x,y) \geq 0$ for all $x,y$, it follows
by Slepian's lemma that per choice of ${\bj}$, the probability
$q^b_{n,\bs,\bj}$ is maximal when $\ss_m=1$ for all $m$.
In case $n$ is odd, note that
$f^b_{n+1}(x) = \sigma_{n}(x) \hat{f}^b_n(x) + b_n x^n$ and
$$
|x|^n \leq 2 n^{-\delta/2} \sigma_{n}(x)  \quad \quad \forall \ x \in V
$$  
Consequently, for all ${\bf j}$ and $\bs$,
\beq{temp}
q^b_{n,\bs,{\bf j}} \leq P(|b_n| \geq n^{\delta/4}) +
P\Bigl( \min_{m=1}^{k+4} \inf_{x \in L_m} \ss_m \hat{f}_n^b (x) >
 - 2 n^{-\delta/8}\Bigr)\;.
\eeq
With $n-1$ even, continuing as before, we see that the right-most term in
(\ref{temp})
is maximal when $\ss_m=1$ for all $m$. 
In conclusion, it suffices to consider 
$$
q^b_{n,\bj} = P\Bigl( \inf_{x \in V'} \hat{f}^b_{n+1} (x) > -2n^{-\delta/8}
\Bigr)\,,
$$
for $n$ even. Applying the arguments of 
\req{4.2}, \req{4.4} and \req{ofer-5.5}
with 
$\gamma_n(x) \equiv -n^{-\delta/8}$ on
the subset $V'$ of $V$, we find that
\beq{6.2}
q^b_{n,\bj} \leq 
e^{-n^{\delta/4}/2}
+ \prod_{i=1}^4 P\Bigl(\sup_{t \in \cT_i} Y_t \leq \ep_{T'}\Bigr) \,,
\eeq
where $Y_t$ is the stationary Gaussian process of Lemma \ref{lem-2.2} and
for $i=1,\ldots,4$ the set $\cT_i \subset [\delta T,(1-\delta)T]$ is the
image of
$V' \cap \cI_i$ under the transformation $t=-\log(1-\theta_i(x))$.
Since $\tau \mapsto R_y(\tau)$ is monotonically decreasing on
$[0,\infty)$, it follows by Slepian's
lemma that $P(\sup_{t \in \cT_i} Y_t \leq \ep_{T'})$ is maximal
per fixed size of $\cT_i$ when the latter set is an interval, that is,
when the $J_{j_l}$ are all at one end of $\cI_i$ for each $i$
(easiest to see this by considering first $J_{j_1}$ only, then $J_{j_2}$ etc.).
In this case each interval $\cT_i$ has at least the length $(1-2\delta)T-k$,
so the upper bound of \req{6.1} follows from \req{2.6} and
\req{6.1new}--\req{6.2}.
\qed

\section{Proof of Theorem \ref{theo-2}}\label{sec5}

\setcounter{equation}{0}

\noindent
In view of the upper bound $q_{n,k} \leq n^{-b+o(1)}$ of
Section \ref{sec6}, it
suffices to show that $p_{n,k} \geq n^{-b-o(1)}$, in order
to complete the proof
of Theorem \ref{theo-2}.

Since $a_i$ are of zero mean and positive variance, the support 
of their law must intersect both $(0,\infty)$ and $(-\infty,0)$. 
Consequently, there exist $\beta<0<\alpha$ such that
$P(|a_i -\alpha| < \epsilon)>0$ and $P(|a_i - \beta| < \epsilon)>0$
for all $\epsilon>0$. Replacing $\{a_i\}$ by $\{-a_i\}$ does not
affect the number of zeros of $f_{n+1}(x)$. Hence,  
we may and shall assume without loss of generality that 
$|\alpha| \ge |\beta|$.
Let $s \ge 4$ be an even integer such that $\alpha+(s-1)\beta<0$.
Define
$$
	Q(x) =\beta x^{s-1} + \sum_{i=0}^{s-2} \alpha x^i,   \quad \quad
	R(x) =\alpha + \sum_{i=1}^{s-1} \beta x^i,
$$
and note that
\beq{QRlemma}
Q(x)>0  \quad \forall |x| \leq 1, \quad \quad R(1)<0<R(-1)\,.
\eeq

\subsection{Proof for $k$ and $n$ even}
\label{keven}

Suppose that $k$ and $n$ are even.
After $k,s,\alpha,\beta$ are fixed, we shall 
choose $\delta>0$ sufficiently small, then a large enough integer 
$r=r(\delta)$, followed by a small enough positive
 $\epsilon=\epsilon(\delta,r)$,
all of which are independent of $n$.
Let $r_i$ denote the multiple of $s$ nearest to $r^i$, for $i=1,\ldots,k$ and
$\rho_n:=\max\{5/p_n,(\log n)^{-1/2}\}$ for $p_n \uparrow \infty$ such that 
$E|a_i|^{p_n} \leq n$ 
(these choices are slightly different from the ones
made in Section \ref{sec7}). 
Let $m=m_n$ be the multiple of $s$ nearest to
$2 r_k \rho_n \log n/|\log (1-\delta)|$
and define the polynomials
\beaa
B(x) = \sum_{i=0}^{m-1} b_i x^i &:=&
	(1+x^s+x^{2s}+\dots+x^{r_1-s}) Q(x)
	+ (x^{r_1}+x^{r_1+s}+\dots+x^{r_2-s}) R(x) \\
	& +& (x^{r_2}+x^{r_2+s}+\dots+x^{r_3-s}) Q(x)
	+ (x^{r_3}+x^{r_3+s}+\dots+x^{r_4-s}) R(x) \\
	& +& \dots
	+ (x^{r_k} + x^{r_k+s} + \dots + x^{m-s}) Q(x), \\
C(X) = \sum_{i=0}^{m-1} c_i X^i &:=& (1+X^s+X^{2s}+\dots+X^{m-s}) Q(X).
\eeaa
Each coefficient $b_i$ of $B(x)$ equals either $\alpha$ or $\beta$.
The same holds for each coefficient $c_i$ of $C(X)$.
Let ${\mathcal A}_n$ denote the event that the following hold:
\begin{itemize}
\item[{\bf A1}]\quad $|a_i-b_i|<\epsilon$ for $i=0,1,\dots,m-1$
\item[{\bf A2}]\quad $|a_{n-i}-c_i|<\epsilon$ for $i=0,1,\dots,m-1$
\item[{\bf A3}]\quad $a_m + a_{m+1}x + \dots + a_{n-m} x^{n-2m}
	> n^{-1/4} \sigma_{n-2m+1}(x)$ for all $x \in \reals$.
\item[{\bf A4}]\quad $|a_i|<n^{\rho_n}$ for $i=0,1,\dots,n$.
\end{itemize}

Most of our work shall be to show that the polynomial $B(x)$ has
the required behavior in terms of zeros for $|x| \leq (1-\delta)^{1/r_k}$.
Conditions {\bf A1} and {\bf A4} ensure that
$f(x)$ is close enough to $B(x)$ on this interval
so as to have there exactly $k$ simple zeros.
The condition {\bf A3} precludes additional zeros of $f(x)$
near $\pm 1$. Moreover, with {\bf A2} and the positivity of $C(X)$
for $|X| \leq 1$, we conclude that $f(x)>0$ when $|x| > 1$.

The stated lower bound on $p_{n,k}$
is an immediate consequence of the following two lemmas.
\begin{lemma}
\label{Bnprob}
For any fixed $\delta>0$, $\epsilon>0$ and an integer $r$, 
the probability of the event ${\mathcal A}_n$ is at least $n^{-b-o(1)}$
for even $n \to \infty$.
\end{lemma}
\begin{lemma}
\label{kzeros}
Suppose the even integer $k$ is fixed. There exist small enough $\delta>0$, 
large enough $r=r(\delta)$ and a small enough $\epsilon=\epsilon(\delta,r)$ positive, 
such that for all sufficiently large even $n$,
any polynomial $f(x)=\sum_{i=0}^n a_i x^i$
whose coefficients are in ${\mathcal A}_n$
has exactly $k$ real zeros, each of which is a simple zero.
\end{lemma}

\noindent
{\bf Proof of Lemma \ref{Bnprob}:}
Since all coefficients of $B(x)$ and $C(x)$ belong to $\{\alpha,\beta\}$,  
our choice of $\alpha$ and $\beta$ implies that 
each coefficient condition in {\bf A1} and {\bf A2}
is satisfied with probability
at least $c$ for some $c>0$ depending only on $\epsilon$.
The probability that condition {\bf A3}
holds is $P_{n',\gamma_{n'}}$ of Theorem \ref{the-3}
for $n'=n-2m+1$ odd and
$\gamma_{n'}(x)\equiv n^{-1/4}$ that satisfy the assumptions
of this theorem.  Consequently, condition {\bf A3} holds with probability of
at least $(n')^{-b-o(1)}=n^{-b-o(1)}$.
Since conditions {\bf A1}, {\bf A2} and {\bf A3} are independent,
the probability that all of them hold is at least
\beq{5.1n}
c^{2m} (n')^{-b+o(1)} = n^{-b-o(1)}.
\eeq
(Recall that $\rho_n \to 0$, hence also $m_n/\log n \to 0$.)
By Markov's inequality and the choice of $\rho_n$, the 
probability that condition {\bf A4} fails for a given $i$ is at most $n^{-4}$.
Hence the probability that this condition fails for any
$i$ in the range $0 \le i \le n$ is at most 
$O(n^{-3})$.
Since $b \leq 2$,
imposing condition {\bf A4}
does not affect the $n^{-b-o(1)}$ lower bound of \req{5.1n}.
\qed

\noindent
{\bf Proof of Lemma \ref{kzeros}:}
The proof of the lemma is divided in three steps.

\noindent
{\it Step 1:} For $\delta>0$
sufficiently small, $r> (\log \delta)/(\log (1-\delta))$
sufficiently large, $\epsilon>0$ sufficiently small
and all large even integers $n$, each
polynomial $f(x)$ with coefficients in ${\mathcal A}_n$ has 
exactly $k$ simple zeros in $[0,1]$.

\noindent
{\it Step 2:} Under same conditions on the parameters, $f(x)>0$ on $[-1,0]$.

\noindent
{\it Step 3:} $g(X)=X^n f(X^{-1}) >0$ on $(-1,1)$.

\noindent
{\bf Step 1.} Fixing $f(x)$ as above,
observe that
the zeros of $f(x)$ in $(0,1)$ are the same as those of $F(x):=(1-x^s)f(x)$,
so it suffices to prove that
\begin{itemize}
\item	$F(x)>0$ for $x \in [0,\delta^{1/r_1}]$
\item	$F'(x)<0$ for $x \in [\delta^{1/r_1},(1-\delta)^{1/r_1}]$
\item	$F(x)<0$ for $x \in [(1-\delta)^{1/r_1},\delta^{1/r_2}]$
\item	$F'(x)>0$ for $x \in [\delta^{1/r_2},(1-\delta)^{1/r_2}]$
\item	$F(x)>0$ for $x \in [(1-\delta)^{1/r_2},\delta^{1/r_3}]$
\item	$F'(x)<0$ for $x \in [\delta^{1/r_3},(1-\delta)^{1/r_4}]$ \\
	\vdots
\item	$F'(x)>0$ for $x \in [\delta^{1/r_k},(1-\delta)^{1/r_k}]$
\item	$F(x)>0$ for $x \in [(1-\delta)^{1/r_k},2^{-1/m})$
\item	$f(x)>0$ for $x \in [2^{-1/m},1]$
\end{itemize}
Indeed, the sign changes in $F(x)$ force at least one real zero
in each of the $k$ gaps {\em between} the intervals
on which $F$ is guaranteed positive or negative,
and the monotonicity of $F'(x)$ on these gaps
guarantees that each of them contains exactly one zero 
and that the zero is simple. Note also that 
$m=m_n \to \infty$, so per choice of $\delta>0$ all the 
intervals of $x$ above are nonempty as soon as $r$ is large enough.

Recall that our choice of $m=m_n$ is such that for any $l<\infty$,
\beq{nxbd}
m^l n^\rho (1-\delta)^{m/r_k} \to 0 \;.
\eeq
Consequently, by conditions {\bf A1} and {\bf A4} of ${\mathcal A}_n$,
there exists $c(r,\delta)$ finite, such that for 
all $\epsilon>0$, $n$ large enough and $|x|\leq (1-\delta)^{1/r_k}$,
\beqn{Ferr}
|F(x)-(1-x^s) B(x)|	&\leq& (1-x^s) \left[ \epsilon(1+|x|+\dots+|x|^{m-1})
			+ n^{\rho} (|x|^m+\dots+|x|^n) \right]	\nn \\
	&\le& (1-|x|)^{-1} (\epsilon + n^\rho x^m) \leq c(r,\delta) \epsilon.
\eeqn
Fix $M$ such that $|Q(x)| \le M$ and $|R(x)| \le M$ for all $x \in [0,1]$.
By definition of $B(x)$,
\beq{Ffake}
(1-x^s)B(x) = Q(x)+
	\left[ \sum_{\ell=1}^k (-1)^\ell (Q(x)-R(x)) x^{r_\ell} \right]
		- Q(x) x^m.
\eeq
Suppose $x \in [0,\delta^{1/r_1}]$.
Then, each $x^{r_\ell}$ and $x^m$ is at most $\delta$, so
$$ 
(1-x^s)B(x) \geq Q(x) - (2k+1)M\delta \;.
$$
Therefore, for all $\delta$ sufficiently small,
the positivity of $Q(x)$ on $[0,1]$ (see \req{QRlemma})
implies that $(1-x^s)B(x) \geq \eta$ for some
$\eta>0$ independent of $n$, and all
$x \in [0,\delta^{1/r_1}]$. For $\epsilon>0$ small enough, this in
turn implies the positivity of $F(x)$ on this interval (see \req{Ferr}).

Suppose $x \in [(1-\delta)^{1/r_j},\delta^{1/r_{j+1}}]$
for some $j \in \{1,2,\dots,k-1\}$.
Then, $x^m \leq x^{r_\ell} \leq \delta$ for all $\ell > j$ and
$x^{r_\ell} \in [1-\delta,1]$, for all $\ell \leq j$.
In view of the identity
\beq{Ffake2}
Q(x)+\sum_{\ell=1}^j (-1)^\ell (Q(x)-R(x)) =
Q(x) 1_{j \ is \ even} + R(x) 1_{j \ is \ odd}\,
\eeq
and \req{Ffake}, it follows that for all $x$ as above,
$$
|(1-x^s)B(x) -
[Q(x) 1_{j \ is \ even} + R(x) 1_{j \ is \ odd}]| \leq (2k+1)M\delta.
$$
For $\delta$ small enough, the error $(2k+1)M\delta$ is at most
$\min\{Q(1),-R(1)\}/3$. Once $\delta$ is chosen,
taking $r$ sufficiently large guarantees that
$Q(x) \geq  Q(1)/2$ and $R(x) \leq R(1)/2$ for all
$x \in [(1-\delta)^{1/r_1},1]$. Since $Q(1)$ is positive and $R(1)$ is
negative (see \req{QRlemma}), we conclude that
there exists $\eta>0$ independent of $n$ such that
$(-1)^j (1-x^s)B(x) \geq \eta$ for all $n$ large enough
and all $x \in [(1-\delta)^{1/r_j},\delta^{1/r_{j+1}}]$, $j=1,\ldots,k-1$.
In view of \req{Ferr}, for all $\epsilon>0$ small enough
$(-1)^j F(x)$ is then positive throughout the interval
$x \in [(1-\delta)^{1/r_j},\delta^{1/r_{j+1}}]$, as needed.

Suppose $x \in [(1-\delta)^{1/r_k},2^{-1/m}]$.
Then, $x^{r_\ell} \in [1-\delta,1]$ for all $\ell \le k$ and $x^m \le 1/2$.
With $k$ even, it follows from \req{Ffake} and \req{Ffake2} that
$$
(1-x^s)B(x) \geq \frac{1}{2} Q(x) - 2kM\delta \;.
$$
So, when $\delta>0$ is small enough, then
for some $\eta>0$ independent of $n$, it holds that
$(1-x^s)B(x) \geq \eta$ for all $n$ large enough and all
$x \in [(1-\delta)^{1/r_k},2^{-1/m}]$. Recall that
\beq{mdl}
(1-x^s) x^m(a_m+a_{m+1}x+\dots+a_{n-m} x^{n-2m}) \geq 0
\eeq
by condition {\bf A3}. So, while $F(x)-(1-x^s)B(x)$ is no longer
negligible as in \req{Ferr}, the positivity of
the expression in \req{mdl} results in
$$
F(x)-(1-x^s)B(x) \geq
- (1-x^s) (\epsilon (1+|x|+\dots+|x|^{m-1})+ m (\alpha+\epsilon) |x|^{n-m})
\geq - c(\delta,r) \epsilon
$$
for some finite $c(\delta,r)$, all $\epsilon>0$ and large enough $n$.
(This is because $m_n=o(\log n)$, so
$m_n x^{n-m_n} \to 0$ as $n \to \infty$, uniformly on
$x \in [(1-\delta)^{1/r_k},2^{-1/m}]$).
Consequently, when $\epsilon>0$ is small enough,
the uniform positivity of $(1-x^s) B(x) \geq \eta >0$
results in the positivity of $F(x)$
for $x \in [(1-\delta)^{1/r_k},2^{-1/m}]$.

Suppose $x \in [2^{-1/m},\delta^{1/n}]$. Using the
decomposition $f=f^L+f^M+f^H$ as in \req{LMH}, note that by condition {\bf A1}
\beaa
f^L(x) \geq B(x) - \epsilon \sum_{i=0}^{m-1} |x|^i &\geq&
(x^{r_k}+x^{r_k+s}+\dots+x^{m-s})Q(x) - r_k M - \epsilon m  \\
&\geq& \Big( \frac{m-r_k}{2 s} \Big) Q(x) - r_k M - \epsilon m.
\eeaa
Note that $f^M(x) \geq 0$ by condition {\bf A3} and
$$
|f^H(x)| \leq m (\alpha+\epsilon) x^{n-m} \le m (2\alpha) \delta^{1-m/n}
 \le 2 m \alpha \delta^{1/2}
$$
by condition {\bf A2}. Since $m_n \to \infty$
while $m_n/n \to 0$,
and $Q(x)$ is strictly positive,
we see by combining the above that if $\delta$ and $\epsilon$ are
small enough then for all $n$ large enough the
``main'' term $m_n Q(x)/(2s)$ dominates,
so $f(x) >0$ for all $x \in [2^{-1/m},\delta^{1/n}]$.

Suppose $x \in [\delta^{1/n},1]$. In this case, by condition {\bf A3},
\beq{Fmid}
f^M(x)
\ge x^m n^{-1/4} \sigma_{n-2m+1}(x) \ge
\delta^2 n^{-1/4} \sqrt{n-2m+1} > n^{1/8}
\eeq
as $n \rightarrow \infty$.
Condition {\bf A1} implies that $|f^L(x)| \leq (\alpha+\epsilon) m$,
whereas condition {\bf A2} implies that $|f^H(x)| \leq (\alpha+\epsilon) m$.
Since $m=o(\log n)$, we conclude that $f(x)>0$
for large $n$ and all $x \in [\delta^{1/n},1]$.

\bigskip

We turn to deal with the sign of $F'(x)$ in the gaps
$[\delta^{1/r_j},(1-\delta)^{1/r_j}]$ for $j=1,\dots,k$.
To this end, first note that
$$
F'(x) = \frac{d}{dx}\left[(1-x^s)B(x) \right] +e(x)
$$
where by conditions {\bf A1} and {\bf A4}, there exists $c(\delta,r)$ finite, such that 
for all $\epsilon>0$, $n$ large enough and $x \in [0,(1-\delta)^{1/r_k}]$,
\beqn{dFerr}
|e(x)| &\leq & |-sx^{s-1}| \left[ \epsilon(1+x+\dots+x^{m-1})
			+ n^{\rho} (x^m+x^{m+1}+\dots+x^n) \right]\nn	\\
	&&\nichts+  (1-x^s) \left[ \epsilon(1 + 2x + \dots + (m-1)x^{m-2})
		+ n^{\rho} (mx^{m-1}+(m+1)x^m+\dots+nx^{n-1}) \right]\nn \\
	&\le& s (1-x)^{-1} (\epsilon + n^\rho x^m)
	+ (1-x)^{-2} ( \epsilon + n^\rho m x^{m-1}) \leq c(\delta,r) \epsilon
\eeqn
(see \req{nxbd} and \req{Ferr}). Next, using \req{Ffake}, we obtain
\beq{Ffakederivative}
	\frac{d}{dx}\left[(1-x^s)B(x) \right] = 
	 Q'(x) + \sum_{\ell=1}^k (-1)^\ell \left[ (Q'(x)-R'(x)) x^{r_\ell} + 
(Q(x)-R(x)) r_\ell x^{r_\ell-1} \right]  -  o(1).
\eeq
(The $o(1)$ denotes two terms involving $x^m$, which by 
\req{nxbd} converge to $0$ uniformly on $x \in [0,(1-\delta)^{1/r_k}]$.) 
The sum of the terms involving $Q'(x)$ or $R'(x)$ in 
\req{Ffakederivative} is at most $(2k+2)M'$, where
$M'$ is such that $|Q'(x)| \le M'$ and $|R'(x)| \le M'$ for all $x \in [0,1]$.
Per fixed $\delta>0$, if $r$ is sufficiently large then  
$Q(x)-R(x) \geq \eta$ for some $\eta>0$ and all 
$x \in [\delta^{1/r_1},1]$ (see \req{QRlemma}).
We claim that if $x \in [\delta^{1/r_j},(1-\delta)^{1/r_j}]$
for some $j \in \{1,2,\dots,k\}$ then the term
$h_j:=(-1)^j (Q(x)-R(x)) r_j x^{r_j-1}$
dominates the right hand side of \req{Ffakederivative}
for all $r$ large enough. 
Indeed, $|h_j| \geq \eta \delta r_j$
for all $x \in [\delta^{1/r_j},(1-\delta)^{1/r_j}]$,
whereas for such $x$ we have that 
$|h_\ell| \leq 2 M r_{j-1}$ when $\ell<j$ and
$|h_\ell| \leq 2M r_k (1-\delta)^{(r_{j+1}-1)/r_j}$ when $\ell>j$.
Since $r_\ell = r^\ell (1+o(1))$, combining the above we see that 
for all large enough $r$,
$$
(-1)^j\frac{d}{dx}\left[(1-x^s)B(x) \right] \geq {\eta \over 2} \delta r^j -  
3 M k (r^{j-1}+r^k (1-\delta)^{\sqrt{r}}) - (2k+2) M'  -o(1) 
\geq {\eta \over 3}\delta r^j.
$$
By \req{dFerr} we then get that for small enough $\epsilon>0$,
$(-1)^j F'(x)$ also is positive in the $j$-th gap.

This completes Step 1.

\bigskip
\noindent
{\bf Step 2.}
As before, define $F(x):=(1-x^s)f(x)$.
The proof that $F(x)>0$ on $[-\delta^{1/r_1},0]$
is the same as the proof for $[0,\delta^{1/r_1}]$,
now using the positivity of $Q(x)$ on $[-1,0]$.
For each $j \in \{1,2,\dots,k-1\}$,
the analysis for $[-\delta^{1/r_{j+1}},-(1-\delta)^{1/r_j}]$
is the same as that for $[(1-\delta)^{1/r_j},\delta^{1/r_{j+1}}]$,
the only difference is that $Q(x)$ and $R(x)$ are both positive
near $-1$ (whereas they have opposite signs near $1$),
so the result is that $F(x)>0$ on these intervals,
independent of the parity of $j$.
The analyses for $[-2^{-1/m},-(1-\delta)^{1/r_k}]$
and for $[-1,-2^{-1/m}]$ are the same as for
the symmetric intervals on the positive side.

To complete the proof that $f(x)>0$ on $[-1,0]$,
it remains to show that $F(x)>0$ on each gap
$[-(1-\delta)^{1/r_j},-\delta^{1/r_j}]$ for $j=1,\dots,k$.
By \req{Ferr}, it suffices to show that on such an interval 
$(1-x^s)B(x) \geq \eta$ for some $\eta>0$, independent
of
$\epsilon$ and $n$.
On the $j$-th such interval, $x^m \leq x^{r_\ell} \leq \delta$ for all $\ell > j$,
whereas if $r$ is sufficiently large, then 
$1 \geq x^{r_\ell} \geq \delta^{r_\ell/r_j} \geq (1-\delta)$ for all $\ell < j$.
Hence, it follows from \req{Ffake} and \req{Ffake2} that 
\beq{Ffake3}
	|(1-x^s)B(x) - [t(x) Q(x) + (1-t(x)) R(x)]| \leq (2k-1)M\delta\,,
\eeq
where $t(x)=1-x^{r_j}$ for $j$ even, and $t(x)=x^{r_j}$ otherwise.
Let $\eta=\min\{Q(-1),R(-1)\}/4>0$, and take 
$\delta$ small enough that $(2k-1)M \delta < \eta$.
Since $t(x) \in [0,1]$,
if $r$ is large enough that $\min\{Q(x),R(x)\}>2\eta$ for
all $x \in [-1,-\delta^{1/r_1}]$, then \req{Ffake3} implies
that $(1-x^s) B(x) \geq \eta$ for all $x \in
[-(1-\delta)^{1/r_j},-\delta^{1/r_j}]$, $j=1,\dots,k$. The 
positivity of $F(x)$ for small $\epsilon$ and large $n$
follows (by \req{Ferr}).

\bigskip
\noindent
{\bf Step 3.}
To complete the proof of Lemma~\ref{kzeros}, it suffices to show
that $g(X):=X^n f(X^{-1})$ is positive on $(-1,1)$. For $\epsilon < \alpha$,
conditions {\bf A1}, {\bf A2} and {\bf A3} result in
\beq{ginq}
g(X) \geq C(X) -(\epsilon + 2\alpha X^{n-m}) \sum_{i=0}^{m-1} |X|^i 
               + X^m n^{-1/4} \sigma_{n-m+1}(X) \,,
\eeq
for all $|X| \leq 1$. Since $(1-X^s) C(X)=(1-X^m)Q(X)$, we see that for $n$
large enough and all  $|X| \leq 2^{-1/m}$,
$$
(1-X^s) g(X) \geq (1-X^m) Q(X) - 
(\epsilon + 4 \alpha 2^{-n/m}) \sum_{i=0}^{s-1} |X|^i
\geq \frac{1}{2} Q(X) - 2 s\epsilon 
$$
is positive for $\epsilon < Q(1)/(8s)$. Since $C(X) \geq \frac{m}{2s}Q(X)$
when $|X| \in [2^{-1/m},\delta^{1/n}]$, it follows from \req{ginq} that 
$$
g(X) \geq m \left[\frac{Q(X)}{2s} - \epsilon -  2\alpha \delta^{1-m/n}\right] 
$$ 
is positive for any such $X$, 
provided $\epsilon < Q(1)/(8s)$, $2\alpha \delta^{1/2} < Q(1)/(8s)$
and $n$ is large enough. Finally, for large $n$, if $|X| \in [\delta^{1/n},1]$ then 
$X^m n^{-1/4} \sigma_{n-m+1}(X) \geq n^{1/8}$  (see \req{Fmid}).
Since $m=o(\log n)$,
the positivity of $g(X)$ for such $X$ is a direct consequence of
\req{ginq}.
\qed

\subsection{Proof for $k$ and $n$ odd}

In this section we sketch the modifications
to the argument of the previous section
that are required for the case where $k$ and $n$ are odd.
We will specify an event occurring with probability at least $n^{-b-o(1)}$
that forces $k-1$ simple zeros in $(0,1)$, one simple zero in $(-\infty,-1)$,
and no other real zeros.
Fix positive $\delta$, integer $r$ and $\epsilon>0$, and define 
$r_i$, $\rho=\rho_n$ and $m=m_n$ as in Section \ref{keven}.
Define the polynomials
\beaa
B(x) = \sum_{i=0}^{m-1} b_i x^i &:=&
	(1+x^s+x^{2s}+\dots+x^{r_1-s}) Q(x)
	+ (x^{r_1}+x^{r_1+s}+\dots+x^{r_2-s}) R(x) \\
	& & \nichts + (x^{r_2}+x^{r_2+s}+\dots+x^{r_3-s}) Q(x)
	+ (x^{r_3}+x^{r_3+s}+\dots+x^{r_4-s}) R(x) \\
	& & \nichts + \dots
	+ (x^{r_{k-1}} + x^{r_{k-1}+s} + \dots + x^{m-s}) Q(x), \\
C(X) = \sum_{i=0}^{m} c_i X^i &:=& (1+X^s+X^{2s}+\dots+X^{r_1-s}) Q(X) \\
		& & \nichts + \alpha X^{r_1}
		+ X (X^{r_1}+X^{r_1+s}+\dots+X^{m-s}) Q(X)
\eeaa
the coefficients of which are in $\{\alpha,\beta\}$.
Let ${\mathcal B}_n$ denote the event that the following hold:
\begin{enumerate}
\item[{\bf B1}]	\quad $|a_i-b_i|<\epsilon$ for $i=0,1,\dots,m-1$
\item[{\bf B2}]	\quad $|a_{n-i}-c_i|<\epsilon$ for $i=0,1,\dots,m$
\item[{\bf B3}]	\quad $a_m + a_{m+1}x + \dots + a_{n-m-1} x^{n-2m-1}
	> n^{-1/4} \sigma_{n-2m}(x)$ for all $x \in \reals$
\item[{\bf B4}]	\quad $|a_i|<n^{\rho}$ for $i=0,1,\dots,n$.
\end{enumerate}

Note that the degree of $C(X)$ is one larger than in Section~\ref{keven}.
This ensures that the ``middle polynomial'' in condition 
{\bf B3}
has even degree, so that Theorem~\ref{the-3} applies to it.
Hence, 
similarly to the proof of Lemma \ref{Bnprob}, 
one has that the event ${\mathcal B}_n$ occurs for odd $n$
with probability exceeding $n^{-b-o(1)}$.

For all small enough $\delta>0$, large enough $r$ and small enough
$\epsilon$, 
the argument of the proof of Lemma \ref{kzeros}, using the shape of $B(x)$,
shows that if the coefficients of $f(x)$ are in ${\mathcal B}_n$ 
then $f(x)$ has exactly
$k-1$ zeros in $[0,1]$, all simple, and no zeros in $[-1,0]$.
We next prove that the function $F(X):=(1-X^s) X^n f(1/X)$ 
satisfies
\begin{itemize}
\item	$F(X)>0$ for $X \in (0,1)$
\item	$F(X)>0$ for $X \in [-\delta^{1/r_1},0)$
\item	$F'(X)>0$ for $X \in [-(1-\delta)^{1/r_1},-\delta^{1/r_1}]$
\item	$F(X)<0$ for $X \in [-2^{-1/m},-(1-\delta)^{1/r_1}]$
\item	$F(X)<0$ for $X \in (-1,-2^{-1/m}]$
\end{itemize}
These will imply that $f(x)$ has a simple zero in $(-\infty,-1)$,
and no other zeros with $|x|>1$.
Together with the $k-1$ simple zeros in $[0,1]$,
this will bring the total number of zeros to $k$.

First, a proof analogous to that of~(\ref{Ferr}) shows that
there exists $c(r,\delta)$ finite, such that
for $|X| \le (1-\delta)^{1/r_1}$,
\beqn{newFerr}
	|F(X)-(1-X^s) C(X)| \leq c(r,\delta) \epsilon.
\eeqn
The analogue of~(\ref{Ffake}) is
\beqn{newFfake}
	(1-X^s) C(X) = (1-X^{r_1}+X^{r_1+1}-X^{m+1}) Q(X)
				+ \alpha X^{r_1}(1-X^s).
\eeqn

Suppose $X \in [0,(1-\delta)^{1/r_1}]$.
Then~(\ref{newFfake}) implies
	$$(1-X^s) C(X) \ge (1-X^{r_1}) Q(X) \ge \delta Q(X)
		> c(r,\delta) \epsilon$$
if $\epsilon$ is small enough, so $F(X)>0$ by~(\ref{newFerr}).

Suppose $X \in [-\delta^{1/r_1},0]$.
Then~(\ref{newFfake}) implies
	$$(1-X^s) C(X) \ge (1-3\delta) Q(X)$$
so $F(X)>0$ by~(\ref{newFerr}) assuming suitable $\delta$ and $\epsilon$.

Suppose $X \in [-(1-\delta)^{1/r_1},-\delta^{1/r_1}]$.
The analogues of (\ref{dFerr}) and (\ref{Ffakederivative})
are
	$$\left| F'(X) - \frac{d}{dX}\left[(1-X^s)C(X) \right] \right|
	\le c(\delta,r) \epsilon$$
and,
with $r_1$, $s$ even,
\begin{eqnarray*}
	\frac{d}{dX}\left[(1-X^s)C(X) \right] &=& 
	\left( 1-X^{r_1}+X^{r_1+1} \right) Q'(X)
		+ \left[ -r_1 X^{r_1-1} + (r_1+1) X^{r_1} \right] Q(X) \\
	&& \nichts	+ \alpha r_1 X^{r_1-1} (1-X^s)
		- \alpha s X^{r_1+s-1} - o(1) \\
	&\ge& -3M' + r_1 \delta Q(X) - \alpha r_1 (1-\delta^{s/r_1}) - o(1)\\
	&\ge& (r_1\delta/2) Q(X),
\end{eqnarray*}
in which the last inequality holds for $r$ sufficiently large.
Hence for $\epsilon$ small enough, $F'(X)$ will be positive.

Suppose $X \in [-2^{-1/m},-(1-\delta)^{1/r_1}]$.
Then,
for $r$ sufficiently large, $1-X^{r_1}+X^{r_1+1}-X^{m+1} \leq -(1/2-3\delta)/2$ 
and $\alpha X^{r_1} (1-X^s) = O(s/r_1)$. For $\delta$ small and $r$ large
(\ref{newFfake}) thus implies that
	$$(1-X^s) C(X) \le -(1/2-3\delta)Q(X)/2 + O(s/r_1) \le -Q(X)/8$$
Although~(\ref{newFerr}) is no longer valid, we may apply {\bf B3}
to deduce
\begin{eqnarray*}
	F(X)-(1-X^s) C(X) &\le&
\epsilon s +
	(1-X^s)(a_m X^{n-m} + a_{m+1} X^{n-m-1} + \dots + a_{n-m-1} X^{m+1})\\
	&& \nichts+ m (\alpha+\epsilon) X^{n-m+1} \\
	&\le& \epsilon s + m(\alpha+\epsilon) 2^{-(n-m+1)/m} \leq 2 \epsilon s
\end{eqnarray*}
since $m 2^{-n/m} \to 0$.
Hence $F(X)<0$ 
if we take first $\delta$ small then $r$ large and finally $\epsilon$ small.

Similarly, for $r$ large enough, if $X \in [(1-\delta)^{1/r_1},2^{-1/m}]$, then 
$$(1-X^s) C(X) \ge (1-X^m) X Q(X) \geq Q(X)/3,
$$
and 
$$F(X)-(1-X^s)C(X)  \ge 
	-\epsilon s - m(\alpha+\epsilon) X^{n-m+1} \geq -2 \epsilon s
$$
implying that $F(X)>0$ in this interval.

Both the proof that $F(X)>0$ on $[2^{-1/m},1)$ and the proof
that $F(X)<0$ on $(-1,-2^{-1/m}]$
parallel the proof in Section~\ref{keven}
that $f(x)>0$ for $x \in [2^{-1/m},1)$.
\qed

\section{Proof of Proposition \ref{prop-2}}\label{sec10}

\setcounter{equation}{0}

In view of the upper bound $q_{n,k} \leq n^{-b+o(1)}$ of Section \ref{sec6}, 
it suffices to provide a lower bound
on the probability of the event considered in 
Proposition \ref{prop-2}. To this end, partitioning and shrinking the $U_i$ 
if necessary, we may assume that $m_1=\ldots=m_k=1$, and that the closures 
of the $U_i$ avoid both $1$ and $-1$. Let $\delta \in (0,1/3)$ then be
such that each of
the $U_i$ is contained either in $(-1+\delta,1-\delta)$  or its
image under the 
map $\inv(x)=x^{-1}$. Let $r$ be the number of
$U_i$ of the former type and $s=k-r$ the number of those of the latter
type. Let $\SSS=(-\eta/2,\eta/2)$ for $\eta>0$ as in the statement of the 
proposition. Fix the polynomials $B(x)=\sum_{i=0}^r b_i x^i \in \SSS [x]$
and $C(X)=\sum_{i=0}^s c_i X^i \in \SSS [X]$ with coefficients in $\SSS$,
such that $B(x)$ has $r$
real zeros, one in each of the $U_i$ that are contained in $(-1,1)$ whereas
$C(X)$ has $s$ real zeros, one in $\inv(U_i)$
for each $U_i$ contained in
$(-\infty,-1) \cup (1,\infty)$. Without loss of generality we can set
$b_r>0$ and $c_s>0$. Let $\rho_n=5/p_n$ for $p_n \uparrow \infty$ such that 
$E|a_i|^{p_n} \leq n$ 
(these differ from the quantities defined in
Section \ref{sec7}).
Define the
even integer $m=m_n=2 \lfloor \rho_n \log n/|\log (1-\delta)| \rfloor$
depending on $n$.
For fixed $\epsilon \in (0,\eta/11)$,
consider the event ${\mathcal C}_n$ that all of the following are satisfied:
\begin{itemize}
\item[{\bf C1}] \quad  $|a_i-b_i|<\epsilon$ for $0 \le i \le r$, \quad
   $|a_{r+i}- 9\epsilon 1_{i \; {\rm even}}\,| < \epsilon$ for $0<i<m$.
\item[{\bf C2}] \quad  $|a_{n-i}-c_i|<\epsilon$ for $0 \leq i \leq s, \quad$
   $|a_{n-s-i}- 9\epsilon 1_{i \; {\rm even}}\,| < \epsilon$ for $0<i<m$.
\item[{\bf C3}] \quad 
$a_{m+r}+a_{m+r+1}x+\dots+a_{n-s-m} x^{n-k-2m} >
0$
for all $x \in \reals$.
\item[{\bf C4}] \quad  $|a_i| < n^{\rho_n}$ for $0 \le i \le n$.
\end{itemize}

Proposition \ref{prop-2} is an immediate consequence of
the following two lemmas.
\begin{lemma}
\label{lem-5.1}
For any fixed $B(x)$, $C(X)$ with coefficients in $\SSS$ 
and positive $\epsilon < \eta/11$,
the probability of the event ${\mathcal C}_n$ is at least $n^{-b-o(1)}$.
\end{lemma}
\begin{lemma}
\label{lem-5.2}
For fixed $B(x)$ and $C(X)$, if $\epsilon>0$ is sufficiently small 
and $n$ sufficiently large,
then any polynomial $f(x)=\sum_{i=0}^n a_i x^i$ satisfying
the conditions of ${\mathcal C}_n$ has exactly $k$ real zeros, one in each of
the $U_i$ intervals.
\end{lemma}
{\bf Proof of Lemma \ref{lem-5.1}:} Note that
$P(a \in G)>0$ for any open subset $G$ of $(-\eta,\eta)$ (by
our assumption about the support of the law of $a_i$). Hence
each coefficient condition in {\bf C1} and {\bf C2} 
is satisfied with probability
at least $c$ for some $c>0$ depending only on $B(x)$, $C(X)$ and $\epsilon$.
We continue along the lines of the proof of Lemma \ref{Bnprob}
(taking now $n'=n-k-2m+1$
and $\gamma_{n'}=0$).
\qed

\noindent
{\bf Proof of Lemma \ref{lem-5.2}:}
Our choice of $\rho=\rho_n$ and $m=m_n$ guarantees that 
for any $l<\infty$, $m^ln^\rho (1-\delta)^m\to_{n\to\infty} 0$.
Consequently,
by {\bf C1} and {\bf C4},
for some $\kappa_0=\kappa_0(\delta)$, 
all $\epsilon>0$, 
$n>n_0$ for some $n_0=n_0(\delta,\ep)$
large enough and
$|x| \leq (1-\delta)$
\beaa
 |f(x)-B(x)|     &\le& 10 \epsilon(1+|x|+\dots+|x|^{m+r-1})
            + n^\rho (|x|^{m+r}+|x|^{m+r+1} + \dots)   \\
   &\le& (10 \epsilon + n^\rho (1-\delta)^{m+r})/\delta \le \kappa_0 \epsilon \;.
\eeaa
Hence if $\epsilon$ is small enough
and $n$ large enough, $f$ must have at least as many zeros
as $B(x)$ within $(-1+\delta,1-\delta)$.
On the other hand, $B^{(r)}(x)$ is a positive constant,
and for $x \in (-1+\delta,1-\delta)$,
\beaa
    |f^{(r)}(x)-B^{(r)}(x)|
    &\le& 10 \epsilon \sum_{i=r}^{m+r-1} i^r |x|^{i-r}
            + n^\rho \sum_{i=m+r}^{\infty} i^r |x|^{i-r}  \\
    &\le& 10 \epsilon \sum_{i=r}^\infty i^r (1-\delta)^{(i-r)}
            + n^\rho \sum_{i=m+r}^{\infty}
                    i^r (1-\delta)^{(i-r)},
\eeaa
which again can be made arbitrarily small by shrinking $\epsilon$.
So, we can and shall assume $f^{(r)}(x)>0$ in $(-1+\delta,1-\delta)$.
By Rolle's Theorem, this bounds the number of real zeros in $(-1+\delta,1-\delta)$
by $r$, so $f(x)$ has exactly $r$ zeros in $(-1+\delta,1-\delta)$.
Moreover, taking $\epsilon>0$ such that $|B(x)|>\kappa_0 \epsilon$ for all
$|x| \leq (1-\delta)$, $x \notin U_i$, $i=1,\ldots,k$, implies that the constant
sign of $f(x)$ between each adjacent pair of intervals $U_i$ that are contained
in $(-1+\delta,1-\delta)$ is the same as the sign of $B(x)$ there. Hence $f(x)$
has exactly one zero in each of the $r$ intervals $U_i$ contained in $(-1,1)$.
Similar
 arguments (using {\bf C2} and {\bf C4}) show that for
 some $\kappa_1=\kappa_1(\delta)$ and all $|X| < (1-\delta)$,
$$
    |X^n f(X^{-1})- C(X)|   \le \kappa_1 \epsilon,
$$
with the $s$-th derivative of the polynomial $X^n f(X^{-1})$ made positive
throughout $|X| < (1-\delta)$ by shrinking $\epsilon$.
Recall that $C(X)$ has exactly one zero in each of the intervals
$\inv(U_i)$ for the $U_i$ contained in $(-\infty,-1)\cup (1,\infty)$.
Thus,
for small enough $\epsilon$, the same property holds for the
$s$ zeros of $X^n f(X^{-1})$ within $|X| < (1-\delta)$.

It thus remains to show that
$x^r f(x)>0$ for $(1-\delta) \leq |x| \leq (1-\delta)^{-1}$.
Since $2r+m$ is an even integer, we have by condition {\bf C3}
that for all $x \in \reals$,
\beq{bdf3}
x^r \left(\sum_{i=m+r}^{n-s-m} a_i x^i \right)
= x^{2r+m} \left( \sum_{i=0}^{n-k-2m} a_{m+r+i} x^i \right)
\geq 0 \,.
\eeq
It is easy to check that
$$
h_m(x):= 8 x^2 \sum_{j=0}^{(m-4)/2} x^{2j} -
|x| \sum_{j=0}^{(m-2)/2} x^{2j} \geq 0 \,,
$$
for all even $m \geq 4$ and $2/3 \leq |x| \leq 3/2$. Consequently, for
$(1-\delta) \leq |x| \leq (1-\delta)^{-1}$ and $\delta < 1/3$,
by condition {\bf C1},
\beq{bdf1}
x^r \left(\sum_{i=r+1}^{r+m-1} a_i x^i\right) 
	= x^{2r} \left(\sum_{i=1}^{m-1} a_{r+i} x^i\right)
\geq \epsilon x^{2r} h_m (x) \geq 0,,
\eeq
whereas for $r+(n-s)=2r+n-k$ an even integer, by condition {\bf C2},
\beq{bdf2}
x^r \left(\sum_{i=n-s-m+1}^{n-s-1} a_i x^i \right) =
x^{r+(n-s)} \left(\sum_{i=1}^{m-1} a_{n-s-i} x^{-i} \right)\geq
\epsilon x^{2r+n-k} h_m(x^{-1}) \geq 0.
\eeq
Next note that for sufficiently small $\epsilon>0$,
the polynomial $\sum_{i=0}^r a_i x^i$ has a
positive leading coefficient and no zeros for $|x| \geq (1-\delta)$,
so $x^r (\sum_{i=0}^r a_i x^i)>0$ for all $|x| \ge (1-\delta)$.
Similarly,
$x^{-(n-s)}\sum_{i=n-s}^n a_i x^{i}$ is then a polynomial
with positive constant coefficient
and no zeros for $|x| \leq (1-\delta)^{-1}$. With $r+(n-s)$
an even integer, it
follows that $x^{r} (\sum_{i=n-s}^n a_i x^{i}) \geq 0$ for
$|x| \leq (1-\delta)^{-1}$. In view of \req{bdf3}--\req{bdf2},
we find that $x^r f(x)>0$ for $(1-\delta) \leq |x| \leq (1-\delta)^{-1}$.
\qed

\end{document}